\documentclass[journal,twosides,web]{IEEEtran}

\usepackage[latin1]{inputenc}
\usepackage{ae}
\usepackage{amssymb}
\usepackage{graphicx,epsfig,psfrag}
\usepackage{amsmath}
\usepackage{mathtools}
\usepackage{rotate}
\usepackage{color}
\usepackage{amsfonts}
\usepackage{setspace}
\usepackage{fancyhdr,enumerate}

\usepackage{enumitem}
\usepackage{pifont}
\usepackage{marginnote}

\makeatletter
\let\MYcaption\@makecaption
\makeatother
\usepackage[font=footnotesize]{subcaption}
\makeatletter
\let\@makecaption\MYcaption
\makeatother

\usepackage{etoolbox}
\makeatletter
\@ifundefined{color@begingroup}%
  {\let\color@begingroup\relax
   \let\color@endgroup\relax}{}%
\def\fix@ieeecolor@hbox#1{%
  \hbox{\color@begingroup#1\color@endgroup}}
\patchcmd\@makecaption{\hbox}{\fix@ieeecolor@hbox}{}{\FAILED}
\patchcmd\@makecaption{\hbox}{\fix@ieeecolor@hbox}{}{\FAILED}

\usepackage{multirow}
\usepackage[lowtilde]{url}
\usepackage[hidelinks]{hyperref}
\usepackage{booktabs}
\usepackage{makecell}
\usepackage[referable]{threeparttablex}

\newtheorem{problem}{Problem}
\newtheorem{definition}{Definition}

\newtheorem{example}{Example}
\newtheorem{remark}{Remark}
\newtheorem{theorem}{Theorem}
\newtheorem{lemma}{Lemma}

\newtheorem{assumption}{Assumption}
\newtheorem{policy}{Policy}

\newcommand {\be}{\begin{equation}}
\newcommand {\ee}{\end{equation}}
\def\qed{{$\hfill\blacksquare$}}

\def\@IEEEheaderstyle{\normalfont\sffamily\tiny}
\begin{document}
	
\title{An agent--based decentralized threshold policy \\finding the constrained shortest paths}

\author{
Francesca Rosset, Raffaele Pesenti, and Franco Blanchini
\thanks{This work has been submitted to the IEEE for possible publication. Copyright may be transferred without notice, after which this version may no longer be accessible.}
\thanks{F. Rosset and F. Blanchini are with the Dipartimento di Matematica e Informatica, Universit\`a di Udine, 33100 Udine, Italy, (email: blanchini@uniud.it, rosset.francesca@spes.uniud.it)}
\thanks{R. Pesenti is with the Dipartimento di Management, Universit\`a Ca' Foscari Venezia, Italy, (email: pesenti@unive.it)} 
}
	
\maketitle

\begin{abstract}
We consider a problem where autonomous agents enter a dynamic and unknown environment described by a network of weighted arcs. 
These agents move within the network from node to node according to a decentralized policy using only local information, with the goal of finding a path to an unknown sink node to leave the network.
This policy makes each agent move to some adjacent node or stop at the current node. The transition along an arc is allowed or denied based on 
a threshold mechanism that takes into account the number of agents already accumulated in the arc's end nodes and the arc's weight.
We show that this policy ensures path-length optimality in the sense that, in a finite time, 
all new agents entering the network reach the closer sinks by the shortest paths. 
Our approach is later extended to support constraints on the paths that agents can follow. 
\end{abstract}
	
	

\section{Introduction}\label{sec:introduction}
\IEEEPARstart{I}{n} this work, we consider an environment described by a directed network containing \emph{source} nodes and \emph{sink} nodes. 
Memoryless agents (\emph{tokens}) are generated in source nodes and move, using local information, to find a minimum cost path to sink nodes to leave the network. The topology of the network is unknown to the tokens. The cost of each arc, a finite integer, is also unknown until they reach the arc's tail node. 
 
We will show that, in the long run, the tokens can find their minimum cost paths by implementing a simple decentralized \emph{threshold policy}, provided that the network does not contain non-positive cost cycles and that the tokens can stop and be buffered at the nodes.
	
\begin{figure}[htpb]
			\centering
			\includegraphics[width=7cm]{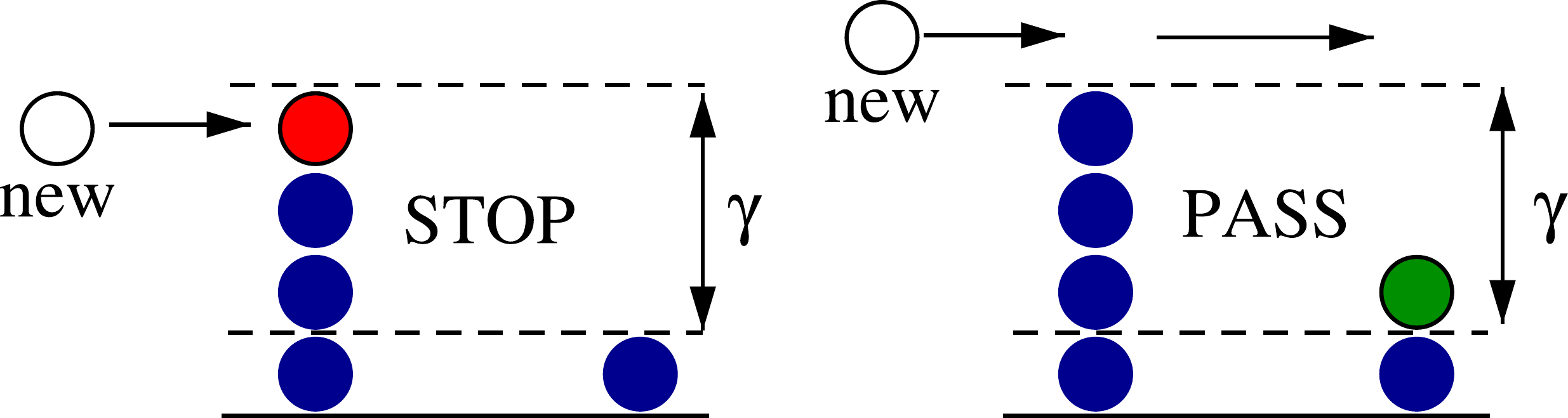}
			\caption{The mechanism with threshold  value $\gamma=3$. The transition between
				two nodes is controlled by a threshold mechanism: no transition is possible
				if the difference between the number of tokens in the buffers of the two nodes is not above a threshold equal to the arc cost (left). 
				A token can pass only when its arrival breaks such a threshold (right).}
			\label{fig:threshold}
		\end{figure}

This policy assumes that when a token arrives at a node~$i$, it finds out which arcs originate from the node, their \emph{cost} $\gamma$, and the \emph{state} $x$ of their extreme nodes (i.e., the number of tokens buffered there). It then immediately decides whether to stop there or to traverse an arc $(i,j)$ to reach the adjacent node~$j$.
The token makes its decision only on the basis of the arc cost~$\gamma_{ij}$ and the tokens stopped in the buffers of the extreme nodes, i.e. the states $x_i$ and $x_j$ of nodes $i$ and $j$.  It applies the {threshold policy} represented in Fig.~\ref{fig:threshold}. First, it increments by one the value of the buffer of node $i$. 
Then, it traverses the arc $(i,j)$ only if this new value exceeds the value of the buffer of the adjacent node $j$ of more than the \emph{threshold} given by the cost~$\gamma_{ij}$ of the arc. The token stops in node $i$ otherwise. 
Formally, if $$x_i + 1 -x_j>\gamma_{ij}$$ the \textit{transition} of the token from~$i$ to~$j$ is \textit{permitted}, otherwise the transition is \textit{denied}.
We say that a token is \emph{above-threshold} in node~$i$ if a transition is permitted for some node~$j$ adjacent to node~$i$, it is \emph{under-threshold} otherwise. In this latter case the token stops in~$i$ and the node state is incremented by one: $x_i \leftarrow x_i +1$.

We assume that each arc $(i,j)$ of the network can also have a finite integer secondary cost $\sigma_{ij}$ independent of $\gamma_{ij}$, which can impose a constraint on the tokens' movements. 
The arcs' costs may represent, e.g., the arcs' lengths, traversal times, or the consumed/recovered energy needed to traverse the arcs. Then, they can be used to force the tokens, e.g., to take a route that minimizes the consumed energy, while not exceeding a given traversal time.
We will prove that our threshold policy can be adapted to deal also with the secondary costs.

Our threshold policy is optimal ``in the long run''. In an initial transient phase, some tokens may not be able to reach the sinks and accumulate at the nodes.  However, we will prove that after the transient, any new token:
	\begin{itemize}
		\item leaves the network and does not increase the number of tokens stopped at each node.
		We say that the network has reached a \emph{steady-state} in terms of tokens accumulated in the nodes;
        \item flows through a shortest path from its source to a sink.
	\end{itemize}
Note that the possibility that the token travels indefinitely in the network is excluded, since we assume the absence of cycles with a non-positive sum of costs.  The accumulation of infinitely many tokens in a node is also excluded by the fact that sooner or later some threshold will be broken.
	
Although the routes found by the tokens turn out to be the least cost (\emph{shortest}) paths, 
the aim of this work is not to propose a new method for solving the well-known shortest path problem.
Its main contribution is to show that the application of simple local rules results in a self-organising emergent behaviour that leads to the discovery of the shortest paths \cite{Bla19,Fig16},
even when an additional path constraint makes the shortest path problem hard to solve.
As additional contributions, we will show that: i) our policy supports negative cost arcs; ii) the steady state of the system does not depend on the initial state of the nodes; iii) the shortest path property is guaranteed even with dynamic networks (e.g., under arc removal/failure or insertion), as long as the network remains connected and its configuration is maintained for a sufficiently long time interval; iv) a bound on the number of tokens needed to complete the initial transient phase can be established. 

In summary, we propose a local policy that yields an exact optimal solution for both the 
shortest path problem and its constrained version.
As our policy is decentralized, it has the following advantages: 
 i) it is fault-tolerant; 
ii) it supports large unknown networks; 
 iii) it does not require notifying the entire network when a change occurs in there. 
 
\subsection{Related works}
The {\em shortest path problem} (SPP) under known data has a long history~\cite{Eva17}, many variants and solution algorithms. It can be solved in polynomial time if there are no negative cycles~\cite{Mad17}. Arc costs are usually assumed to be integers: this is advantageous in terms of space efficiency, speed, and stability of solution algorithms~\cite{Hav19}.
	
The standard SPP algorithms consider static, fully known networks.
In contrast, we deal with an unknown, possibly dynamic environment. 
As mentioned above, our agents can only observe the network at a local level, and the structure of the network may change over time.

The literature proposes different solutions for dealing with unknown \cite{Pap91,Fel04} and dynamic networks \cite{Sun21}. 
We focus on decentralized, agent-oriented methods that use only local information. 
Some optimal decentralized methods based on consensus have been proposed in \cite{Tah06} and \cite{Zha17}. A method based on reinforced random walks is studied in \cite{Fig16}. 
In addition, many metaheuristics, which by definition cannot guarantee optimality,
have been proposed: ant colonies \cite{Dor06,Lis15}, river formation dynamics \cite{Rab17,Red17}, amoebas \cite{XZha17}, 
genetic algorithms \cite{Zhu14}, particle swarm optimization \cite{Moh08}. 
	
Finally, many shortest-path routing protocols exist for packet-switched networks \cite{Wed06}, ad-hoc and mesh networks \cite{Alo12}, and wireless sensor networks (WSN) \cite{Yu14,Ala20}. 
These methods support multipath routing, scalable performance, self-organizing behavior, locality of interaction and network failure detection and backup \cite{Ala20}. The network architecture we consider is indeed like the one in WSN; however, we neglect aspects such as link congestion and throughput.
	
Our approach is inspired by the decentralized threshold policy presented in \cite{Bla19}. It introduces a problem of decentralized control of a continuous flow in a single--source--single--sink network and shows that the flow eventually concentrates along the shortest path.
This mechanism can be used to model phenomena such as lightning~\cite{Bla20}, or electronic circuits of Zener diodes or nonlinear resistors~\cite{Bu99}. 
	
The \textit{(resource) constrained shortest path problem} (CSPP) \cite{Jok66} is an extension of the SPP where an additional constraint on the sum of the secondary costs of the traversed arcs is imposed, making some paths of the network infeasible.  
The standard algorithms for solving the CSPP consider static, known networks.
Since the CSPP is NP-hard \cite{Pug13}, both exact solution methods and heuristics have been proposed. The former are typically based on combinatorial approaches, see the survey in \cite{Fes15}, which when applied to practical problems can yield optimal solutions even in short time \cite{Ahm21}. The latter include various bio-inspired approaches, such as the hybrid particle swarm optimization-variable neighbourhood search algorithm in \cite{Mar17}, but also problem specific heuristics as in \cite{Cab20} or, more recently, the very first learning approaches as in \cite{Yin23}.

The rest of the paper is organized as follows. In Section \ref{sec:def} we define the problem we consider. 
In Section~\ref{sec:MTSDynSys} we present the proposed policy for the simpler unconstrained problem. 
In Section~\ref{sec:constrPaths} we extend these results to the case of the constrained problem.  
In Section~\ref{sec:examples} we report some illustrative examples. 
In Section~\ref{sec:conclusion} we draw some conclusions.
The proofs are in the Appendix.

\section{Definitions and Problem statement}\label{sec:def}
	
A directed network $\mathcal{G}=(\mathcal{N},\mathcal{A})$ is given, where $\mathcal{N} =\{1,2,\ldots,n\}$ is a set of $n$ nodes and $\mathcal{A}\subseteq \mathcal{N}^2$ is a set of $m$ directed arcs. 
Each arc $(i,j)\in \mathcal{A}$ has a finite \emph{cost} $\gamma_{ij}\in\mathbb{Z}$ and a finite \emph{secondary cost} $\sigma_{ij}\in\mathbb{Z}$. 
We call the system \emph{unconstrained} when the secondary cost $\sigma_{ij}$ is not used, and \emph{constrained} otherwise. 
The network $\mathcal{G}$ has a subset of \emph{source nodes} $\mathcal{S} \subseteq \mathcal{N}$ where the tokens are injected from the external environment. The tokens move along the arcs $\mathcal{G}$ according to their direction to reach another subset of \emph{sink nodes} $\mathcal{T} \subseteq \mathcal{N}$, with $\mathcal{T} \cap \mathcal{S} = \emptyset$, where they are ejected into the external environment. 
A \emph{path} $p = \{h_s \in \mathcal{N}, s=1,2\dots r\}$  
connecting a \emph{starting node} $i$ and a \emph{end node} $j$ on $\mathcal{G}$ is an ordered sequence of $r$ non-repeating nodes $i=h_1,h_2 \dots h_r=j$, such that $(h_s,h_{s+1}) \in \mathcal{A}$. We have a \emph{walk} if some of these nodes are repeated; a \emph{circuit} if only the first and last nodes of the sequence coincide.
A path is said \emph{outgoing}, if its last node is a sink.\\ 
Given a path $p$, we call \emph{length} of $p$ the sum of its arc costs:
\begin{equation} \label{eq:lenght}
    L(p) \doteq \sum_{h_k,h_{k+1}\in p}\gamma_{h_kh_{k+1}};
\end{equation}
and \emph{secondary cost} of $p$ the sum of its secondary arc costs:
    \begin{equation} \label{eq:cost}
        C(p) \doteq \sum_{h_k,h_{k+1}\in p}~~\sigma_{h_kh_{k+1}}.
    \end{equation}
In a constrained system, a path $p$ is said \emph{feasible} if
\begin{equation}
    C(p) \leq C_{max},
    \label{eq:path_constr}
\end{equation}
where $C_{max}\in\mathbb{N}$ is the maximum secondary cost allowed for a path, and the partial subpaths $p_i = \{h_1, ..., h_i\}$ from $h_1$ to the intermediate node $h_i$, with $i<r$ are feasible, too (this second condition is automatically satisfied only if $\sigma_{ij}\geq0$).
A path is said \emph{non-positive} if its length is non-positive.
These definitions hold for walks and circuits, too.

\begin{example}
    Consider the network $\mathcal{G}$ in Fig.~\ref{fig:incidenceMatrixNetwork}. 
    The path $p=\{1,2,4,5\}$ from source node $1$ to sink node $5$ has length $L(p)=1+3+0=4$ and secondary cost $C(p)=1+1+0=2$. 
    This is the shortest constrained path when $C_{max}=2$. Instead,
    $p=\{1,2,3,4,5\}$ is the shortest unconstrained path with $L(p)=3$ (and $C(p)=3$).
    \begin{figure}[htpb]
        \centering
        \includegraphics[scale=0.8]{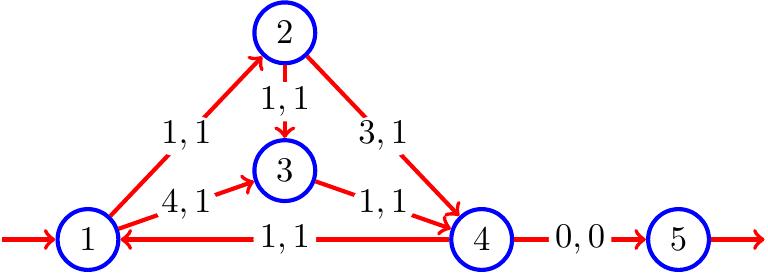}
        \caption{A simple network $\mathcal{G}$. $\gamma_{ij}$ and $\sigma_{ij}$ are shown for each arc $(i,j)$.}
        \label{fig:incidenceMatrixNetwork}
    \end{figure}
\end{example}
	
Each node buffer can store tokens according to a last-in-first-out (LIFO) policy (but our results can be easily generalized to other policies).
In the unconstrained system, we define the {\em state} $x_i(t_k)$ {\em of a node} $i \in \mathcal{N}$ at time $t_k \geq 0$ as the number of tokens present in the node $i$ at time $t_k$.
Hereinafter, for short of notation we write $x(k)$ for $k=0, 1,\ldots$ to mean $x(t_k)$ at time $t_k$ for $k = 0,1,\ldots$.
Note that $x_i(k)=0$ for all $i\in \mathcal{T}$ for all $k \geq 0$ by definition of sink node.
In the constrained system, each node~$i$ buffers tokens in different buckets according to their \emph{current secondary cost} $c$, that is the sum of secondary costs that they paid along the paths they used to reach $i$; 
so the {\em state of a node} $i$ is now a vector $x_i=[x_i^0, x_i^1,\dots,x_i^{C_{max}}]$ with $C_{max} + 1$ components.
Each component $x_i^c(k)$, for $c  = 0,\ldots,C_{max}$,  is the number of tokens with current secondary cost $c$ present in the node $i$ at time $k$. Again, $x_i^c(k)=0$ for all sinks $i\in \mathcal{T}$, for all $0\leq c\leq C_{max}$, for all $k\geq0$.
Note that, in the unconstrained system, all paths are actually feasible, so tokens do not need to know $c$.

We say that a policy is \emph{decentralized} if, at each time, a token makes its routing decision on the basis of the knowledge of the state of the node it currently occupies, the state of the neighboring nodes, and the costs of the arcs incident in the node, but with no memory of the nodes that it has visited so far.
In the constrained case, a token also takes into account the secondary costs of incident arcs and it has memory the secondary cost $c$ it has paid so far. 
Note that the routing decision can be thought to be performed by the ``node agent''
without assuming any computing capability of the token agent.\\
Because nodes buffer tokens according to a LIFO policy, a token arriving at a node has two choices: either remain there or move immediately to a next node. 
We call \emph{destination node} of a token the last node, possibly a sink node in $\mathcal{S}$, reached by a token at which it stops.
	
Accordingly, hereinafter, the subset of nodes that can be reached by a token in $i$ traversing at most one arc is denoted by
		$\mathcal{N}_i = \{j \in \mathcal{N}: (i,j) \in \mathcal{A}\} \cup \{i\}$, for all $i \in \mathcal{N}$.
	
	Let us now introduce the problem of interest.
	\begin{problem}\label{pr:shortestPathConstr}
		Define a decentralized policy that, 
		in the long run, makes each token injected in the nodes $\mathcal{S}$ reach nodes in $\mathcal{T}$, 
		along a feasible path of minimum length and of secondary cost not greater than $C_{max}$.
	\end{problem}
We will first consider its  unconstrained version:
	\begin{problem}\label{pr:shortestPath} Solve Problem \ref{pr:shortestPathConstr} disregarding the secondary costs $\sigma_{ij}$.
	\end{problem}

Throughout this paper, the following assumptions hold.
\begin{assumption}\label{ass:nonNegCircuits}
    Network $\mathcal{G}$ is weakly connected and has:
    \begin{enumerate}
        \item for each source node $i \in \mathcal{S}$, at least one feasible path to a sink node that, in the case of constrained systems, also has a secondary cost less than or equal to~$C_{max}$;
        \item no non-positive circuits with respect the costs $\gamma_{ij}$ and no negative circuits with respect the costs $\sigma_{ij}$; 
        \item the costs $\gamma_{ij}$ and $\sigma_{ij}$ of its arcs $(i,j) \in \mathcal{A}$ upper bounded by a value $\bar{\gamma}$ and $\bar \sigma$, respectively;
        \item  no negative paths between any pair of sinks in $\mathcal{T}$.
    \end{enumerate}
\end{assumption}
\begin{assumption}\label{ass:asincronous}
    Token movements are asynchronous and fast: once a token is injected in the network, it reaches the destination node before any other token is injected.
\end{assumption}
\begin{assumption}\label{ass:memory}
    Tokens are memoryless of the traveled path. Only in  a constrained system, they mantain memory of the cumulative secondary cost $c$ of their traveled paths.
\end{assumption}

In Assumption~\ref{ass:nonNegCircuits}, Condition 1) ensures the existence of a solution to Problem~\ref{pr:shortestPath};
Condition 2) prevents the tokens from looping indefinitely in the network
and ensures the existence of a constrained shortest path (and not a walk) between any source and the sinks;
Condition 4) ensures that a token reaching a sink is expelled immediately, as there are no shorter outgoing paths to other sinks. 
Assumption~\ref{ass:asincronous} implies that there is no coordination on the tokens' movements and that the time taken by tokens to move in the network is of orders of magnitude smaller than the time intervals between successive injections of tokens in the network. 
In practice, it
can be relaxed by not allowing multiple tokens traveling simultaneously in adjacent nodes.
We will revise this assumption in Section~\ref{sec:MTSDynSys} by introducing a multiple-time-scale dynamical system.

	
	\section{The unconstrained problem}\label{sec:MTSDynSys}
	In this section, we present the threshold policy for the unconstrained  Problem~\ref{pr:shortestPath}. 
	To this end, we call \emph{state of the network}~$\mathcal{G}$ the integer vector 
	$x=[x_1,x_2,\dots,x_{|\mathcal{N}|}]^\top\in \mathbb{Z}^{|\mathcal{N}|}$ whose components are the state of the network nodes and we introduce the following definition.
		\begin{definition}[Admissibility]
			The state $x$ of a network $\mathcal{G}$ is {admissible} if $x_i-x_j \leq \gamma_{ij}$ for all $(i,j)\in\mathcal{A}$.
		\end{definition}
        By the next 
		lemma, $x$ is admissible if for any two nodes
		$i$ and $j$ there is no directed path between them with cost smaller than 
		the difference of the state components $x_i-x_j$.
		
		\begin{lemma}\label{lem:adm_path}
			The state $x$ of a network $\mathcal{G}$ is {admissible} if and only if 
			\begin{equation}\label{eq:path}
				x_i-x_j \leq L(p),
			\end{equation}
			for all $i, j \in \mathcal{N}$, for all paths $p$ connecting $i$ to $j$. \qed
		\end{lemma}
    Hereafter, until otherwise stated we assume that the network initial state is admissible.
    
		Our system is multiple-time-scale dynamical. 	
		A coarser {\em slow dynamic} time-scale $t_0, t_1, \ldots, t_k, \ldots$ 
			marks the times in which tokens are injected in the source nodes (i.e., at $t_k$ the $k$-{th} token enters in $\mathcal{G}$).
		  A finer {\em fast dynamic} time-scale $t_k =t_k^0<t_k^1<...<t_k^{N_k}<t_{k+1}$  marks the times of the $N_k$ 
		  \emph{elementary transitions} that a generic token $k$ performs from $t_k$ until it stops or reaches a sink. 
		  We define \emph{elementary transition} the movement of a token from a node $i$ to a node $j$ along the arc $(i,j)$. 
		  We define \emph{network transition} the overall elementary transitions occurring between two consecutive time instants $t_k$ and $t_{k+1}$ of the slow dynamics.
		\begin{assumption}[Slow-fast dynamic system]\label{ass:MTS}
			Network~$\mathcal{G}$ has two kinds of  time dynamics:
			\begin{enumerate}
				\item \emph{Fast dynamics}. Each token moving in the network~$\mathcal{G}$ takes a negligible time 
				to complete its path, as long as elementary transitions are possible. 
				\item \emph{Slow dynamics}. After a token $k$ is injected at time $t_k$, a new token $k+1$ can be injected 
				at time $t_{k+1}$ with the interval $t_{k+1}-t_k$ long enough to ensure that the 
				token $k$ has reached its destination node.
			\end{enumerate}		
    \end{assumption}
		\begin{remark}\label{ass:onetokenmov}
			Assumption~\ref{ass:MTS} implies that at each time at most a token is moving in the network. When a token enters in the network~$\mathcal{G}$, it sees that all the previously injected tokens either already left~$\mathcal{G}$ or are 
			stopped at the nodes. It sees that no other token can move except itself as it is the only one that can be above-threshold in a node.
			Therefore, to study how the tokens accumulates in the nodes or find their paths to the sink nodes,
			we can  consider only the network state of~$\mathcal{G}$  at the times when a token has reach its destination and a new token can enter in~$\mathcal{G}$.
		\end{remark}
		\begin{remark}
			Assumptions \ref{ass:MTS} can be partially relaxed  (although it is out of the scope of this paper) by just assuming negligible arc traversal times, but allowing simultaneous movements of different tokens and new token injections before a moving token reaches its destination.
			In this case, we  should consider how the network state varies considering all the elementary transitions.
			Then, if two tokens have to move from the same node, giving precedence to one or to the other might
			make a difference in the resulting network state. 
			In any case, our optimality results still hold even if their proofs become much more cumbersome.
		\end{remark}

Let $e_i\in \mathbb{Z}^{|\mathcal{N}|}$ be the $i$th canonical basis vector, whose components are $1$ at the $i$th position, and $0$ otherwise.
We now introduce in simple words the transition rule of our system.
We will formalize it immediately after.

\noindent
\textbf{Transition rule} {\em Slow dynamics perspective} -
A token injected in a source node when the network state is equal to~$x$
undergoes a sequence of elementary  transitions until it remains above-threshold in the nodes that it visits. 
Eventually, either it reaches a node~$j$ where it is under-threshold or it leaves the network. The network state is updated as $x\leftarrow x+e_j$ in the former case, it does not change, $x\leftarrow x$, in the latter case. 

\noindent
{\em Fast dynamics perspective} - Any elementary transition along arc $(i,j) \in \mathcal{A}$ is described by the control:

$$
		u_{ij} = \begin{cases}
			1 & \text{if~the transition occurs,} \\
			0 & \text{otherwise,}
		\end{cases}
$$
		and the corresponding nodes' state variation are $x_j \leftarrow x_j+u_{ij} \text{ and }  x_i \leftarrow  x_i-u_{ij}$.
		Similarly, the injection of a token in a source node $i$ is described by the input:
  
		$$
		\nu_{i} = \begin{cases}
			1 & \text{if~the injection occurs,} \\
			0 & \text{otherwise,}
		\end{cases}
		$$
		and results in $x_i \leftarrow x_i+\nu_i$.

The above transition rule can be implemented by a \emph{threshold policy} formalized as follows.
		\begin{policy}[Threshold policy for the unconstrained system]\label{pol:threshold} Consider the following conditions:
  \vspace{-3mm}
		\begin{subequations} \label{eq:threshold}
			\begin{align}
				\intertext{a) at most a token can enter any node $i \in \mathcal{N}$, i.e.,}
				\label{eq:threshold_2}&\sum_{j: i \in \mathcal{N}_j} u_{ji} + \nu_i \leq 1;\\
				\intertext{b) at most a token can leave any node $i \in \mathcal{N}$, i.e.,}
				\label{eq:threshold_3}&\sum_{j \in \mathcal{N}_i} u_{ij} \leq 1;\\
				\intertext{c) a transition along an arc $(i,j) \in \mathcal{A}$ may occur only if the arrival of a token in $i$ makes the difference 
					between the number of tokens present in $i$ and in $j$ exceed the value of the arc cost $\gamma_{ij}$ (threshold), i.e., only if}
				\label{eq:threshold_1}&x_i + \sum_{l: i \in \mathcal{N}_l}u_{li} + \nu_i - x_j > \gamma_{ij},~\forall(i,j) \in \mathcal{A}.
				\end{align}
		\end{subequations}
		Then, the control specifying the policy for  any arc $(i,j)$ is 
		\begin{equation}\label{eq:threshold_u}
	    	 u_{ij}= 
			\begin{cases}
				1 & \text{if~\eqref{eq:threshold_2}~and~\eqref{eq:threshold_3}~and~\eqref{eq:threshold_1}~hold,}\\
				0 & \text{otherwise.}
			\end{cases}
		\end{equation}
		\end{policy}

		Policy \ref{pol:threshold} enjoys the following properties:
		\begin{itemize}
				\item it is decentralized, as each transition is decided on the basis of local information;
				\item it imposes that all the walks followed by tokens are paths.
				More formally, the subnetwork of~$\mathcal{G}$ induced by the arcs~$(i,j)$ such that $u_{ij}=1$ 
                        is acyclic, due to~\eqref{eq:threshold_2}, \eqref{eq:threshold_1}, \eqref{eq:threshold_u}, and Assumption~\ref{ass:nonNegCircuits}.2 
                        (see proof of next Theorem~\ref{th:proc_defined})
                        and, moreover, is composed of not-intersecting paths, due to~\eqref{eq:threshold_2} and~\eqref{eq:threshold_3}.
			\end{itemize}
		  
	 As we observe the system from a slow dynamics perspective, we
  see all the controls $u_{ij}$ applied ``instantaneously''. Hence, we have the following state equation for each node $i$ of~$\mathcal{G}$ such that  $i \not \in \mathcal{T}$: 
		\begin{IEEEeqnarray}{lcl}
	\label{eq:stateEqNode}
	x_i(k+1) &~=~& x_i(k) - \sum_{\mathclap{\substack{j\in \mathcal{N}_i}}} u_{ij}(k) + \sum_{\mathclap{\substack{j:i\in \mathcal{N}_j}}} u_{ji}(k) + v_i(k),~~
		\end{IEEEeqnarray}
		 and $x_i(k) = 0$ for each node $i\in  \mathcal{T}$, with $k= 0,1, \ldots$.
		
		Dynamics~\eqref{eq:stateEqNode} can be compactly defined by the equation
		\begin{equation}\label{eq:state}
			x(k+1) =  x(k) + Bu(k) + v(k),
		\end{equation}
	where $B$ is the incidence matrix of network~$\mathcal{G}$ deprived of the rows associated with the sink nodes, $u(k)=[u_{ij}:(i,j)\in\mathcal{A}]^\top$, $v(k)=[v_1,\dots,v_n]^\top$.
		Equations of this type have been used for flow networks and manufacturing
		\cite{BPRU}.
		
		Hereinafter, we denote by $(\mathcal{G},f)$ a system made of a network~$\mathcal{G}$ with dynamics given by Policy \ref{pol:threshold} and state equation~\eqref{eq:stateEqNode}, satisfying Assumptions~\ref{ass:nonNegCircuits}-\ref{ass:memory}.                                  
		We denote by:
		\begin{itemize}
			\item $V(x(k)) = \sum_{i \in \mathcal{N}} x_i(k)$ the total number of tokens in the network at time $k$. 
			\item $\mathcal{O} \subset \mathbb{Z}^{|\mathcal{N}|}$ the set of the admissible states.
			\item $\partial \mathcal{O} \subseteq \mathcal{O}$ the \emph{rest set}, as subset of $\mathcal{O}$
			of all fixed points of \eqref{eq:stateEqNode}.
		\end{itemize}
		
		\begin{definition}
			We say that $\partial \mathcal{O}$ is a \emph{rest set} if for $x \in \partial\mathcal{O}$ 
			there exists a node $i\in\mathcal{S}$ such that adding a token in it ($v(k) = e_i$)
			leaves the state unchanged, i.e.,  $x(k)=x(k+1)=x$. 
			We say that $\bar x \in \partial \mathcal{O}$ is a \emph{global rest state}
			if adding a token in any node in $\mathcal{S}$ leaves the state unchanged.
		\end{definition}
		We remark that in general there are multiple global rest states $\bar x$.
			However, with a little abuse, we consider a global rest state as a steady-state since the number of tokens stopped in the network nodes cannot change.
	
		The next lemma and theorem respectively state that the set $\mathcal{O}$ of admissible states is bounded and
		that, for any state $x(k) \in \mathcal{O}$, the well posedness property holds, i.e.,
		a new well-defined state $x(k+1) \in \mathcal{O}$ is reached.
		\begin{lemma}[Boundedness of the admissibility states]\label{lem:Vbounded}
			In a system~$(\mathcal{G},f)$  with $x(k)\in\mathcal{O}$, then $V(x(k)) \leq \bar{\gamma} |\mathcal{N}|(|\mathcal{N}|-1)/2, \forall k>0$.\qed
		\end{lemma}
		\begin{theorem}[Well posedness and positive invariance]\label{th:proc_defined}
		    A system~$(\mathcal{G},f)$ initialized with $x(0) \in \mathcal{O}$ is well posed
		    and the set~$\mathcal{O}$ of the admissible states 
			is positively invariant and finite.
			Specifically, for each $k \geq 0$, either $x(k+1) = x(k)$ or $x(k+1) = x(k) + e_j$ for some $j \in \mathcal{N}$.\qed
		\end{theorem}

	\begin{definition}
		We say the injection of tokens in a source is \emph{persistent} if the number of injected tokens  is unbounded.
	\end{definition}

		The next theorem shows that, for any initial state $x(0) \in \mathcal{O}$, a global rest state is reached after the insertion of a finite number of tokens.
		
		\begin{theorem}\label{th:proc_defined_2} 
			A system~$(\mathcal{G},f)$ initialized with $x(0) \in \mathcal{O}$ 
			always reaches a state in~$\partial \mathcal{O}$ in finite time $\bar k$ if at each time $k = 0, \ldots, \bar{k}-1$ 
			a new token is injected in the network (under persistent injection). 
			In particular, the system reaches a global rest state if a sufficiently high number of
			tokens are inserted in {\em all} the source nodes. For a single-source network, 
			the maximum number of new tokens that needs to be inserted to reach a global rest state is
            $$
            \bar{\gamma}|\mathcal{N}|(|\mathcal{N}|-1)/2 - V(x(0)),
            $$
            where $V(x(0))$ is the initial number of tokens present in the nodes' buffers, i.e. $\bar k \leq \bar{\gamma}|\mathcal{N}|(|\mathcal{N}|-1)/2 - V(x(0))$.
            \qed
		\end{theorem}
	
	An immediate consequence of the above theorem is that 
		the time (measured in tokens) to reach a {global} rest state  is
		exponential in the size of the problem input.
		Specifically, the time is 
		pseudo-polynomial, since it is polynomial in the numeric value of the input parameter~$\bar \gamma$,
		which in turn means that the time is exponential in the number of bits of the input necessary to express the numeric value of~$\bar \gamma$.
		\bigskip
	
		So far, we have considered a network with an admissible initial state, $x \in \mathcal{O}$.
		Now, assume that the network state  is  non-admissible, $x \not\in \mathcal{O}$. 
			We can still apply the same Policy \ref{pol:threshold} under Assumption \ref{ass:MTS}.
			However, compared to the previous case, multiple tokens may be above-threshold and could move. 
			To allow only one moving token at a time, a priority criterion is adopted. 
		We introduce some intermediate time instants $t_{k,r}$ of the slow dynamic in-between $t_k$ and $t_{k+1}$, such that $t_k<t_{k,1}<t_{k,2}<...<t_{k,r}<...<t_{k+1}$, in which no token is injected, but the moving token is one already present in the network, which takes priority.
		When a token starts moving from node $i$ and stops in $j$ the state is updated as $x\leftarrow x-e_i+e_j$, while if it leaves the network as $x\leftarrow x-e_i$.
		The next theorem ensures that an admissible state in $\mathcal{O}$ is eventually reached in finite time.

		\begin{theorem}[Attractiveness of $\mathcal{O}$]\label{th:proc_defined_2a} 
			A system~$(\mathcal{G},f)$ initialized with $x(0) \not\in \mathcal{O}$ 
			always reaches a state in~$\mathcal{O}$ in finite time 
			if no new token is injected into the network.\qed
		\end{theorem}
	
		\begin{remark}[Negative $\gamma_{ij}$]\label{rem:negstates}
			Policy \ref{pol:threshold} may induce negative components of $x(k)$ 
			in presence of negative arc costs.
			Obviously, such components do not correspond to physical tokens, but to negative 
			``virtual'' ones that need to be referenced in the count when new 
			tokens arrive. For instance, if an arc has cost $\gamma_{ij}=-1$ and $x_i=x_j=0$, then a transition occurs which produces a
			negative virtual token in $i$ and a virtual one in $j$, resulting in $x_i=-1$ and $x_j=1$.
            All the presented theory continues to hold;
   		we only have  to reinterpret the concept       of state as
		  $x_i(k)=r_i(k)-R_i$, 
		  where $r_i(k)$ is the number of physical tokens and $R_i$ is a possibly negative reference of virtual tokens, reached at time $0$ before any token is inserted, which remains unchanged. 
		\end{remark}

		
		\subsection{Main result: optimality of the paths}\label{sec:ShortestPath}

		Recall that systems~$(\mathcal{G},f)$ reach admissible states $x\in\mathcal{O}$ and, then, a global rest state $\bar x\in\partial\mathcal{O}$, in finite time.
		The next theorem ensures that, when $x$ is admissible, an injected token moves to its destination (which is a sink if $x\in\partial\mathcal{O}$, specifically the closest) along a
		shortest path.
		\begin{theorem}[Shortest path]\label{th:shortestPath} 
			A system~$(\mathcal{G},f)$ at time $k$ with $x(k) \in \mathcal{O}$ is given.
			If $v(k)  = e_i$ for some $i  \in \mathcal{N}$, the injected token reaches its destination node $j\in\mathcal{N}$
			following the unconstrained shortest path from $i$ to $j$.
   
   If $j\in\mathcal{T}$ is a sink, it is a closest sink to $i$, i.e., there is no $j' \in \mathcal{T}$ such that $L(p') < L(p)$ being $p$ (respectively $p'$) the shortest path from $i$ to $j$ (resp. $j'$);
				also, $x_i(k) = L(p)$.\qed
		\end{theorem}

  We now introduce the concept of  {\em maximal rest state}, defined as the ({\em unique}) special global rest state $\bar{\bar{x}} = \{\bar{\bar{x}}_i= L(p_i): i \in \mathcal{N}\}$, where 
  $L(p_i)$ is the length of the shortest path from node $i$ to its closest sink, the very $i$ if $i \in \mathcal{T}$. 
		When the network is in this maximal rest state $\bar{\bar{x}}$, a token injected in {\em any} of its node~$i\in\mathcal{N}$  immediately 
			reaches the closer sink along the shortest path of length~$L(p_i)$. We have that: 
			i) $\bar{\bar x}$ depends only on the network structure;  
			ii) for all admissible states $x\in \mathcal{O}$, $x\leq \bar{\bar x}$ component-wise;  
			iii) $\bar{\bar x}$ is a global rest state for any source set $\mathcal{S}$;  
			iv)  for any global rest state $\bar x\in \partial\mathcal{O}$, $\bar x\leq \bar{\bar x}$ component-wise and $\bar x_i=\bar{\bar x}_i$ for all nodes along at least one shortest path connecting each source to the closer sink.
   
		If multiple shortest outgoing paths exist from node $i$ (even to different sinks), all of these have the same length, so $\bar{\bar x}_i$ is still unique. In the maximal rest state, all these multiple shortest paths are active and potentially available for new tokens; conversely, in a generic global rest state, at least one is active, but not necessarily all of them. Moreover, for some nodes, there are necessarily multiple outgoing arcs which are all above-threshold when a token arrives in there. Depending on the policy of choice for the next node to visit in case of multiple arcs above-threshold, some of those shortest paths are never traversed.

\begin{remark} [Dynamic networks]\label{rem:dynnetworks} From Theorems \ref{th:proc_defined_2a} and \ref{th:shortestPath} it follows 
that whenever a network changes its initial configurations and costs in a new one,
then from the (possibly non-admissible) state the network will reach a new global rest state, eventually ensuring optimality.
\end{remark}

\begin{remark}[Non-integer arc costs]\label{rem:rationalcosts}
If there are rational arc costs $\tilde\gamma_{ij}=n_{ij}/d_{ij}\in\mathbb{Q}$, with $n_{ij}\in\mathbb{Z},  d_{ij}\in\mathbb{N}\smallsetminus\{0\}$, our policy can be applied  to the network with integer costs  $\gamma_{ij} = \mu\tilde\gamma_{ij}$, where $\mu = \text{lcm}(\{|d_{ij}|: (i,j)\in\mathcal{A}\})\in\mathbb{N}$ is the least common multiple of the $d_{ij}$: the topology of the shortest paths is the same. 
Still, $\mu$ is a global parameter to be known in advance, and $\gamma_{ij}\gg \tilde\gamma_{ij}$ if $\mu$ is large, so performance degrades (see Theorem \ref{th:proc_defined_2}).
\end{remark}

 \begin{remark}[Enhanced policy]\label{rem:enhpolicy}
 When an injected token reaches a node~$i$ in which no outgoing arc is above-threshold, i.e.,
$
x_i(k)+1-x_j(k)\leq\gamma_{ij}, \forall j\in\mathcal{N}_i,
$
instead of stopping the token, generate a block of $\gamma_{i\hat{j}} + x_{\hat{j}}(k)-x_i(k)$ virtual tokens in~$i$, 
where $\hat{j} = \arg \min_{j\in\mathcal{N}_i}\{\gamma_{ij} + x_{{j}}(k)-x_i(k)\}$, so that 
 $x_i$ becomes
$x_i(k)\rightarrow x_{\hat{j}}(k)+\gamma_{i\hat{j}}$ and the token can move to $\hat{j}$.
The number of tokens to inject to reach a global rest state reduces as   virtual tokens are generated, making this variation independent from $\gamma_{ij}$ and  improving performance significantly.
In a strongly connected network, this enhanced policy  ensures that any injected token always reaches a sink, although possibly along non-optimal longer paths in the initial transient.
\end{remark}

		\section{The constrained problem}\label{sec:constrPaths}
		We now introduce constraints in the routes on network~$\mathcal{G}$ to solve Problem \ref{pr:shortestPathConstr}. 
		For simplicity, we consider $\sigma_{ij}\geq0$ for all arcs: 
		the extension to negative $\sigma_{ij}$ is straightforward. 

        We impose that each node buffers the tokens stopped in there based on their secondary cost $c=C(p)$, with $p$ being the path traveled by such tokens to reach the node. 
		In particular, we consider a multi-component state for each node 
		$i\in\mathcal{N}$: $x_i=[x_i^0,x_i^1,\dots, x_i^{C_{max}}]\in\mathbb{N}^{C_{max}+1}$, 
		where $x_i^c$ is the number of tokens with secondary cost $c$ stopped in  $i$. 
		We redefine the notion of admissible state as follows.
		\begin{definition}[Admissibility, constrained system]
			The state $x$ of a network $\mathcal{G}$
			is admissible if $x_i^c-x_j^{c+\sigma_{ij}}\leq\gamma_{ij}$ for all  $(i,j)\in\mathcal{A}$, for all $c\in\{0,\dots,C_{max}-\sigma_{ij}\}$.
		\end{definition}
		
		To continue using Policy \ref{pol:threshold}, we just redefine vector $u$ which models the fast dynamic. For all $(i,j) \in \mathcal{A}$, the control which describes the possible movement of a token from node $i$ to node $j$ depends now on the constrained  cost $c=C(p)$ of its traveled path $p$ up to $i$, and is denoted by $u_{ij}^c$. Then,  \eqref{eq:threshold_1} is rewritten as follows; a token in node $i$ with secondary cost $c$ can move to node $j$ only if:
		\begin{equation}              
        x_i^{c} + \sum_{l: i \in \mathcal{N}_l}u^{c-\sigma_{li}}_{li} + \nu_i -x_j^{c+\sigma_{ij}} > \gamma_{ij}.
			\label{eq:threshold_4}
		\end{equation}
		
		We say that the token \textit{falls asleep} in $i$ if any tentative to leave the node would result 
		in a violation of \eqref{eq:path_constr}, i.e., if:
		\begin{equation}
			c+\sigma_{ij} > C_{max},~~\forall (i,j)\in \mathcal{A}, j\in\mathcal{N}_i.
			\label{asleep_cond_1}
		\end{equation}
		An asleep token cannot continue moving; thus, it stops definitively in node $i$, increasing $x_i^{c}$. 
A necessary condition to allow the transition of the token along an arc is
		\begin{equation}
			c+\sigma_{ij} \leq C_{max},~~\textnormal{for~some~} j\in\mathcal{N}_i.
			\label{eq:threshold_5}
		\end{equation}
		
		Then,  control \eqref{eq:threshold_u} is rewritten as 
		$$
		u^{c}_{ij} = 
		\begin{cases}
			1 & \text{if ~\eqref{eq:threshold_2} and~\eqref{eq:threshold_3} and~\eqref{eq:threshold_4} and~\eqref{eq:threshold_5} hold, }\\
			0 & \text{otherwise.}
		\end{cases}
		$$
				
		\subsection{The expanded network model}
		
\begin{figure}[tb]
			\centering
			\includegraphics[scale=0.8]{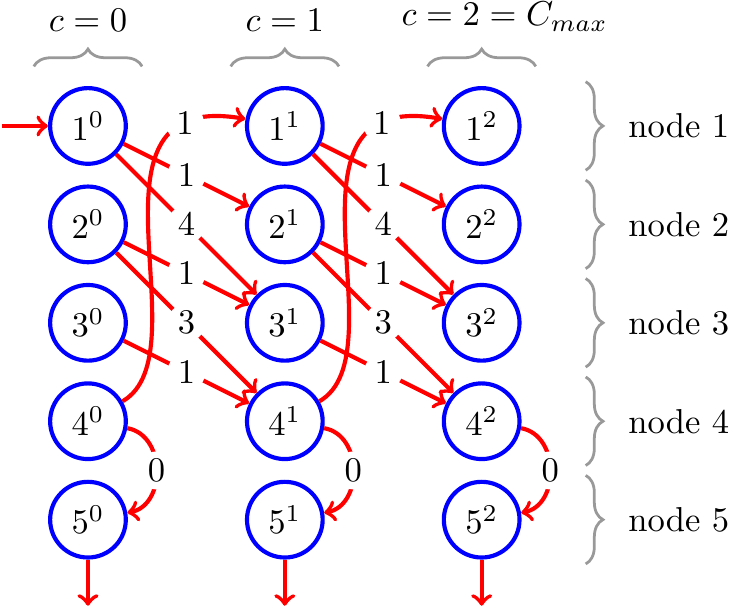} 
			\caption{The expanded network $\mathcal{G}_E$ corresponding to network $\mathcal{G}$ in Fig. \ref{fig:incidenceMatrixNetwork}, when $C_{max}=2$. 
				The value $\gamma_{ij}$ is indicated for each arc $(i,j)$.}
			\label{fig:expanded_network}
		\end{figure}

		To study the constrained system, we model it by an expanded network \cite{Zhu12}, which we denote by $\mathcal{G}_E(\mathcal{N}_E,\mathcal{A}_E)$, 
		composed by $|\mathcal{N}|$ nodes replicated for $C_{max}+1$ steps (each one corresponding to a specific secondary cost $c$), 
		so that $|\mathcal{N}_E|=|\mathcal{N}|(C_{max}+1)$. 
		We denote $i^c$ the $i$-th node of $\mathcal{G}$ replicated at step $c$, for $c\in\{0,\dots,C_{max}\}$. Any arc $(i,j)\in\mathcal{A}$ is replicated, 
		at most, $C_{max}$ times in  $\mathcal{G}_E$, becoming $(i^c,j^{c+\sigma_{ij}})$ for $c\in\{0,\dots,C_{max}-\sigma_{ij}\}$, 
		so that $|\mathcal{A}|\leq|\mathcal{A}_E|\leq|\mathcal{A}|C_{max}$. 
		Note that $|\mathcal{A}_E|=|\mathcal{A}|C_{max}$ when $\sigma_{ij}=1$ for all $(i,j)\in\mathcal{G}$, 
		i.e., we constrain the number of arcs that a token can traverse. 
		Fig. \ref{fig:expanded_network} shows the expanded network of the network in Fig. \ref{fig:incidenceMatrixNetwork}.

		We impose no constraints in the paths of network $\mathcal{G}_E$, so each node has a single-component state. There is a one-to-one correspondence between the state $x_{i^c}$ of any node $i^c\in\mathcal{N}_E$ and $x_i^c$, the $c$-th component of the state $x_i$ of node $i\in\mathcal{N}$. By construction, there is also a one-to-one correspondence  between any unique feasible path or walk $p$ in $\mathcal{G}$ (taking into account its secondary cost $C(p)$), and a unique path $p_E$ in $\mathcal{G}_E$ with the same length $L(p)=L(p_E)$.
		
		If a token falls asleep in node $i\in\mathcal{N}$ of $\mathcal{G}$, this corresponds to reaching a non-sink node $i^{c}\in\mathcal{N}_E$ with no outgoing arcs. For instance, in Fig. \ref{fig:incidenceMatrixNetwork}, a token traveling path $p=\{1,2,3\}$ in $\mathcal{G}$ would fall asleep in node $3$ if $C_{max}=2$. In $\mathcal{G}_E$ (Fig. \ref{fig:expanded_network}), it would reach  node $3^2$, which has no outgoing arcs, through path $p_E=\{1^0,2^1,3^2\}$.
  	Also, some nodes of $\mathcal{G}_E$ are not reachable by the tokens injected in the sources, and can be neglected, as their state components in $\mathcal{G}$ would not be modified. 
			
		\begin{theorem}\label{th:ext_acyclic}
			The  expanded network $\mathcal{G}_E$ of network $\mathcal{G}$:
			\begin{itemize}
					\item is acyclic if and only if there exists no circuit $\varphi$ in $\mathcal{G}$ such that $C(\varphi)=0$;
				  in particular, any walk in $\mathcal{G}$ becomes a valid path in $\mathcal{G}_E$; 
					\item presents a circuit $\varphi_E$ with positive length $L(\varphi_E)>0$ for any circuit $\varphi$ in $\mathcal{G}$ with $C(\varphi)=0$.\qed
			\end{itemize}
		\end{theorem}

		We now prove that a solution for Problem \ref{pr:shortestPath} in $\mathcal{G}_E$ solves Problem \ref{pr:shortestPathConstr} in $\mathcal{G}$.
  Theorem \ref{th:ext_acyclic} ensures no non-positive length circuit in $\mathcal{G}_E$ if Assumption \ref{ass:nonNegCircuits} holds for $\mathcal{G}$, so 
   Assumption \ref{ass:nonNegCircuits} holds for $\mathcal{G}_E$, too. Then,
   in the long run, tokens injected in a source of $\mathcal{G}_E$ ($i^0, i\in \mathcal{S}$) will 
		 reach a closest sink ($j^c,j\in \mathcal{T}$, for some $0\leq c\leq C_{max}$) through the shortest path. 
		
		\begin{theorem}\label{th:constrshortestpaths}
			A system~$(\mathcal{G},f)$ is given.
			A shortest path from a source to a sink in its expanded network $\mathcal{G}_E$ corresponds to a shortest 
			feasible path in $\mathcal{G}$.\qed
		\end{theorem}
  
Note that if $C(\varphi)\geq 0$ does not hold, then the shortest feasible path from a source to a sink in $\mathcal{G}$ might turn out to be a walk (see the proof of the theorem above).

		\begin{remark}
		Applying Policy \ref{pol:threshold} to a constrained system $(\mathcal{G},f)$ is equivalent to consider an unique unconstrained system $(\mathcal{G}_E,f)$. Remarks \ref{rem:negstates}-\ref{rem:enhpolicy} can be adapted to this case.
		\end{remark}

		From the above remark and  Theorem~\ref{th:proc_defined_2},  
		the time to reach a {global} rest state (and  solve the problem) is exponential in the size of the of the problem input.
		Indeed, we must inject a number of tokens that is  polynomial in $\bar \gamma$ in the network $\mathcal{G}_E$ that has a number of nodes polynomial in $C_{max}$.
		So, the time is pseudo-polynomial in the numeric values of both   $\bar \gamma$ and $C_{max}$, i.e., exponential in the number of bits of the input necessary to express them.

		\begin{remark}[Negative $\sigma_{ij}$]
		Negative $\sigma_{ij}$ may make some $x_i^c$ with $c<0$ reachable.  Still, as long as Assumption \ref{ass:nonNegCircuits} holds, there exists $L\leq0$ such that all the  $x_i^c,c<L$ are never reached. The discussion and results  presented above are 
        still valid considering $c\in\{L,...,0,...,C_{max}\}$.
	    \end{remark}


\section{Illustrative examples}\label{sec:examples}
        
        Here, we report two  examples  that apply our policy.
        
        We inject a new token at each time  $k$ of the slow dynamic and  temporize the fast dynamic by taking $|\mathcal{N}|+1$ instants  between $k$ and $k+1$.
        Also, we allow simultaneous transactions involving different nodes, and consider two choice models for the next node to move in when the agent is in a node with multiple outgoing arcs which are above-threshold: i) {\em deterministic}:  the agent  scans the outgoing arcs in a predefined order and chooses the first one that is  above-threshold; ii) {\em stochastic}: the agent chooses one randomly. The initial state is zero, $x_i(0)=0$.
        
        Animations are available in \url{https://users.dimi.uniud.it/~franco.blanchini/examples_mkv.html}.
        
        \subsection{Example 1: a simple network}
       We first study the  network in Fig. \ref{fig:incidenceMatrixNetwork}, with source node $1$ and sink node $5$, using a deterministic choice model.
        
\begin{figure}[t!]
        \centering
        \begin{subfigure}[t]{0.48\linewidth}
            \includegraphics[scale=0.9]{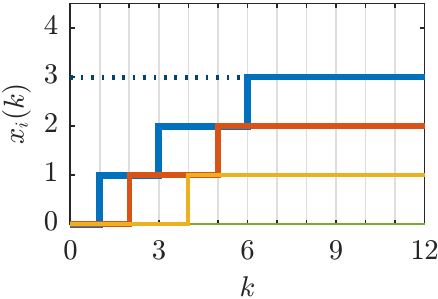}
            \caption{$x_1(k)$ (blue), $x_2(k)$ (red) and $x_3(k)$ (yellow) for the unconstrained system.}
            \label{fig:sim1:time_ev:xunc}
             \vspace*{1mm}
        \end{subfigure}
        \hfill
        \begin{subfigure}[t]{0.48\linewidth}
            \includegraphics[scale=0.9]{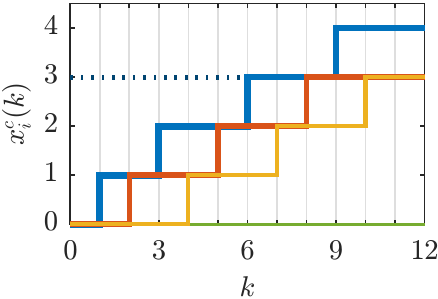}
            \caption{$x_1^0(k)$ (blue), $x_2^1(k)$ (red) and $x_3^2(k)$ (yellow) for the constrained system with $C_{max}=2$.}
            \label{fig:sim1:time_ev:xconcomp}
        \end{subfigure}
        \caption{(\textit{Example 1}). Time evolution of the states of the nodes.  Dotted line: $\bar{\bar x}_1(k)$ for the source of the unconstrained system.} \label{fig:sim1:time_ev} 
        \end{figure}

        Fig. \ref{fig:sim1:time_ev:xunc} shows the time evolution of the states $x_i(k)$ for the unconstrained system. 
        It reaches a global rest state  $\bar x=[3,2,1,0,0]$ after injecting $6$ tokens (which are lost, i.e., stopped in the nodes)  and the shortest path $p=\{1,\allowbreak 2,\allowbreak 3,\allowbreak 4,\allowbreak 5\}$ is eventually discovered and followed by each token injected next. 
        Fig. \ref{fig:sim1:time_ev:xconcomp} shows the evolution of the states $x_i^c(k)$  when constraining paths by $C_{max}=2$. 
        While these reach larger values (non-zero components: $\bar x_1^0=4, \bar x_2^1=\bar x_3^2=3$), the shortest \emph{feasible} path $p=\{1,\allowbreak 2,\allowbreak 4,\allowbreak 5\}$ is  discovered after injecting $10$ tokens (which are lost).
        
        The unconstrained system converges faster, i.e., overall less tokens are required. Indeed, the constrained problem has a larger complexity and can be equivalently seen as an unconstrained problem in the network in Fig. \ref{fig:expanded_network}, which is about $C_{max}$ times larger than the original one in Fig. \ref{fig:incidenceMatrixNetwork}.

        \subsection{Example 2: a large grid  network}\label{sec:sim:2}
\begin{figure}[t!]
        \centering
        \begin{subfigure}[t]{0.48\linewidth}
            \centering
            \includegraphics[width=0.7\linewidth]{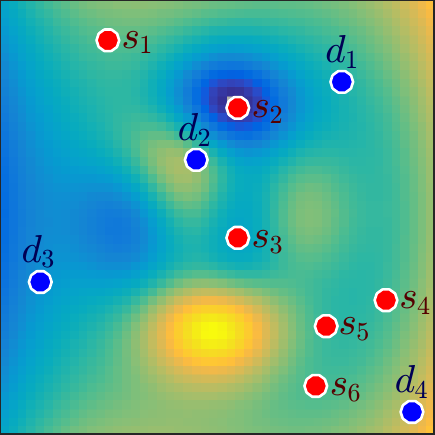}\\
            \vspace*{1mm}
            \includegraphics[scale=0.9]{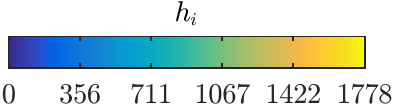}
            \caption{Map of $h_i.$}
            \label{fig:sim2:map:val}
             \vspace*{1mm}
        \end{subfigure}
        \hfill
        \begin{subfigure}[t]{0.48\linewidth}
            \centering
            \includegraphics[width=0.7\linewidth]{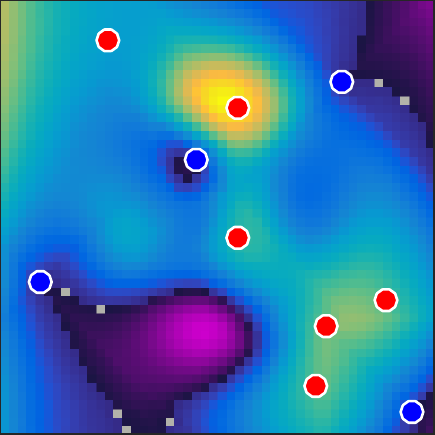}\\
            \vspace*{1mm}
            \includegraphics[scale=0.9]{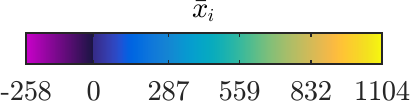}
            \caption{Map of
            $\bar{\bar x}_i$, unconstr. system.}
            \label{fig:sim2:map:barx}
             \vspace*{1mm}
        \end{subfigure}
         
        \begin{subfigure}[t]{0.48\linewidth}
            \centering
            \includegraphics[width=0.7\linewidth]{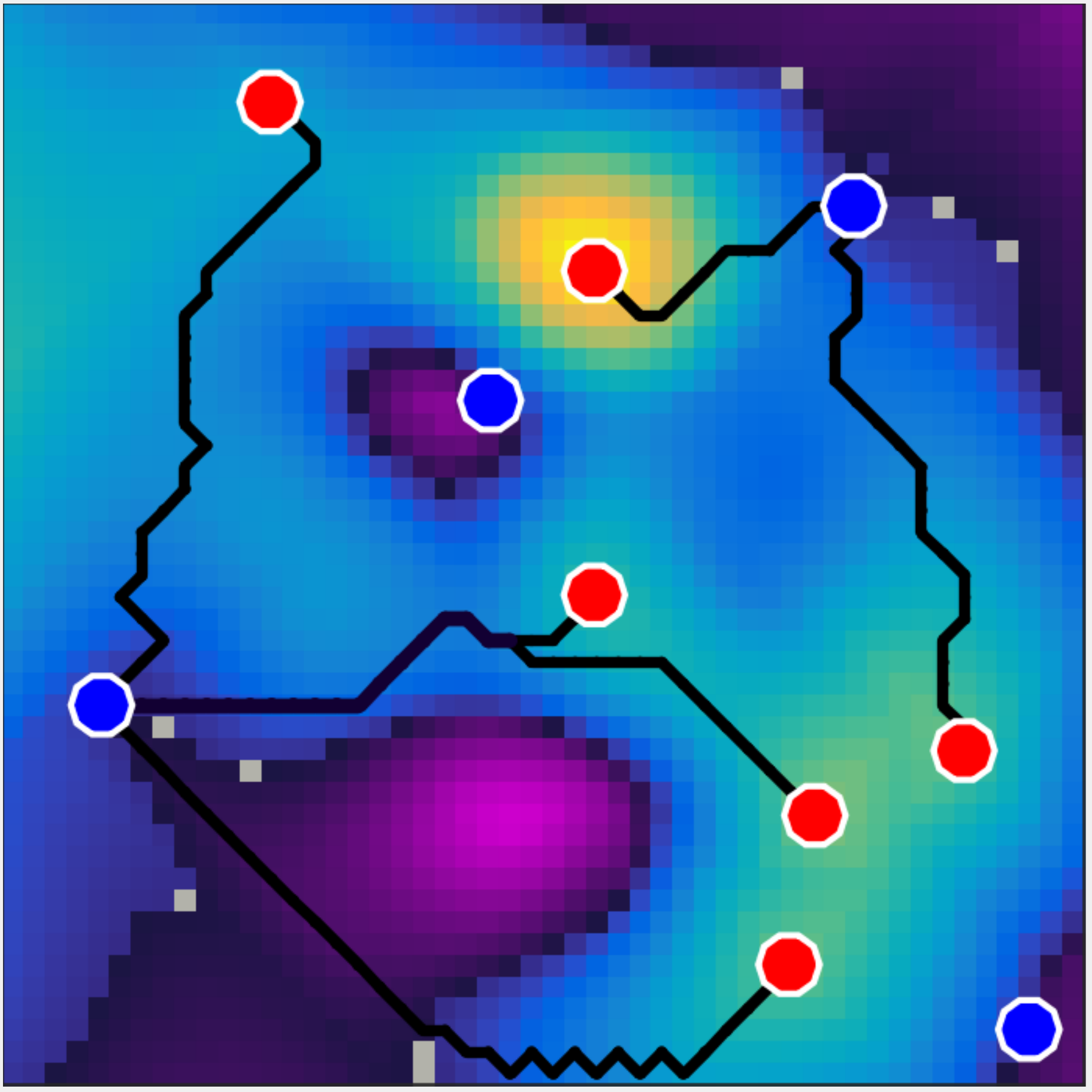}\\
            \vspace*{1mm}
            \includegraphics[scale=0.9]{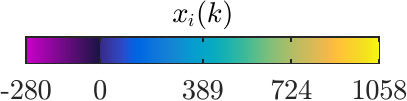}
            \caption{Map of 
            $\bar x_i(k)$, unconstrained system, deterministic choices.}
            \label{fig:sim2:map:xuncdet}
            \vspace*{1mm}
        \end{subfigure}
        \hfill
        \begin{subfigure}[t]{0.48\linewidth}
            \centering
            \includegraphics[width=0.7\linewidth]{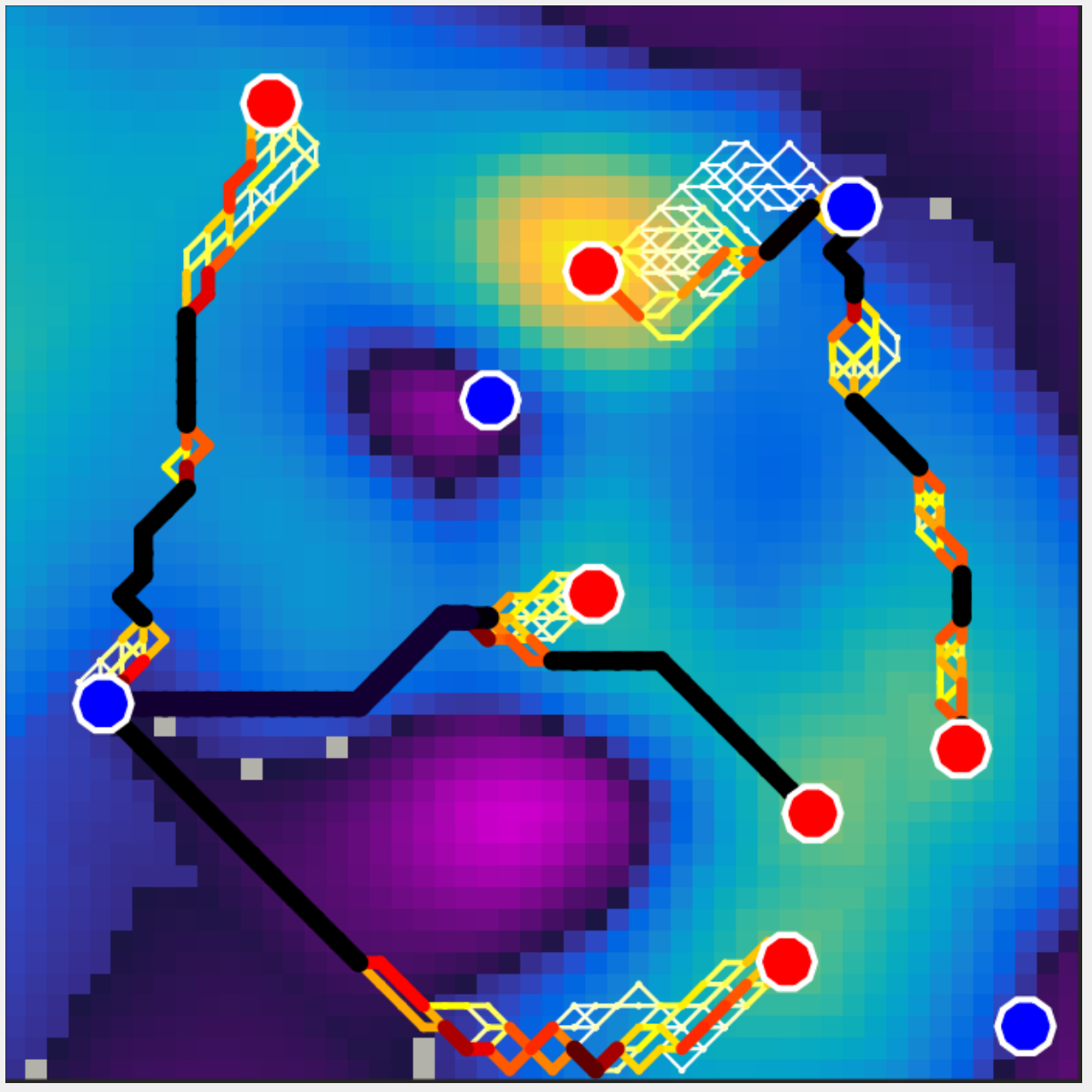}\\
            \vspace*{1mm}
            \includegraphics[scale=0.9]{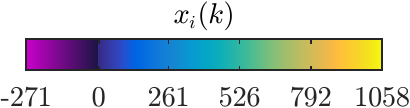}
            \caption{Map of 
            $\bar x_i(k)$, unconstrained system, stochastic choices.}
            \label{fig:sim2:map:xunc}
            \vspace*{1mm}
        \end{subfigure}
         
        \begin{subfigure}[t]{0.48\linewidth}
            \centering
            \includegraphics[width=0.7\linewidth]{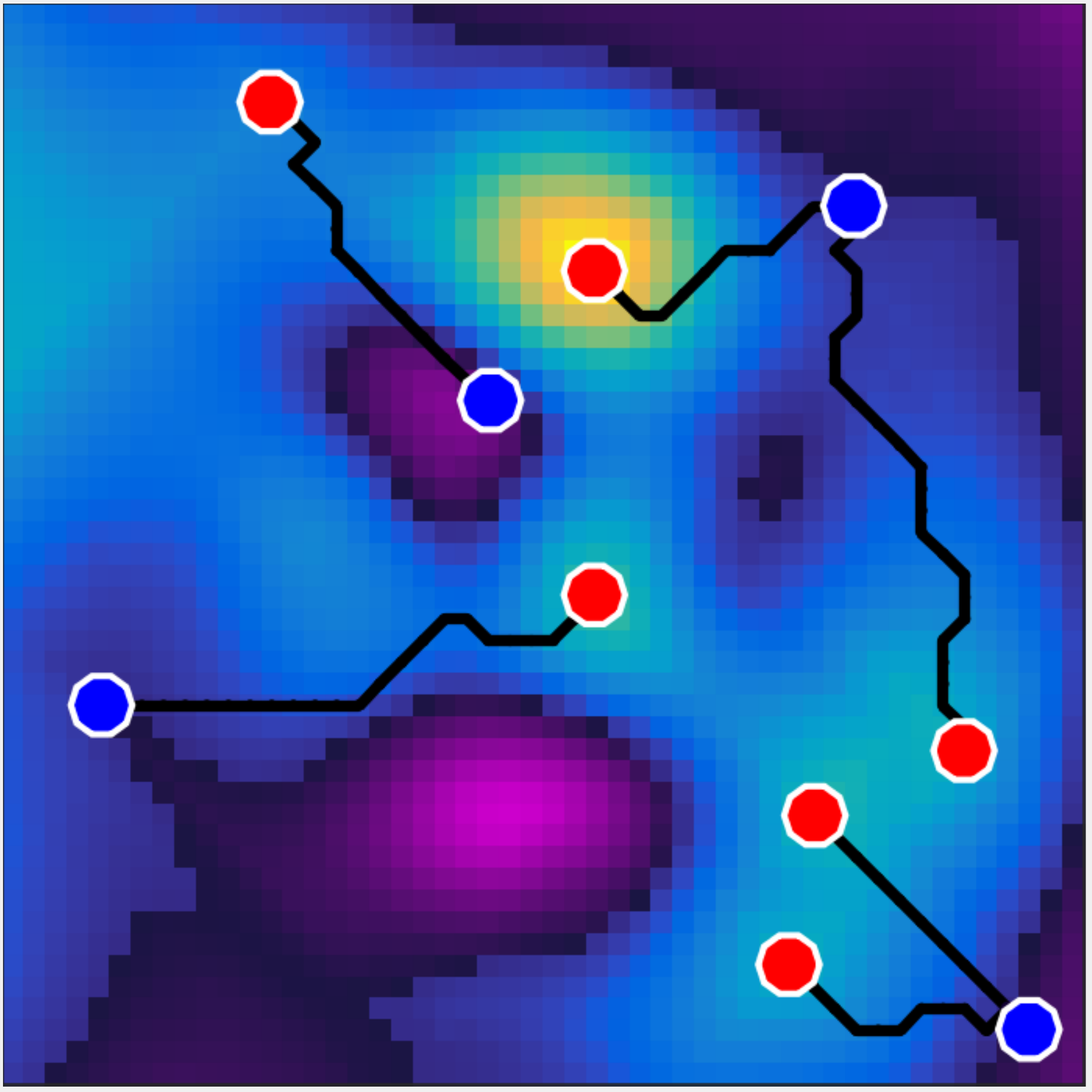}\\
            \vspace*{1mm}
            \includegraphics[scale=0.9]{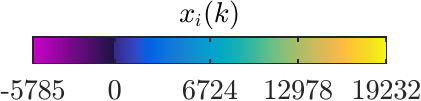}
            \caption{Map of
            $\bar x_i(k)$, constr. system, $C_{max}{=}25$, deterministic choices.}
            \label{fig:sim2:map:xcondet}
        \end{subfigure}
        \hfill
        \begin{subfigure}[t]{0.48\linewidth}
            \centering
            \includegraphics[width=0.7\linewidth]{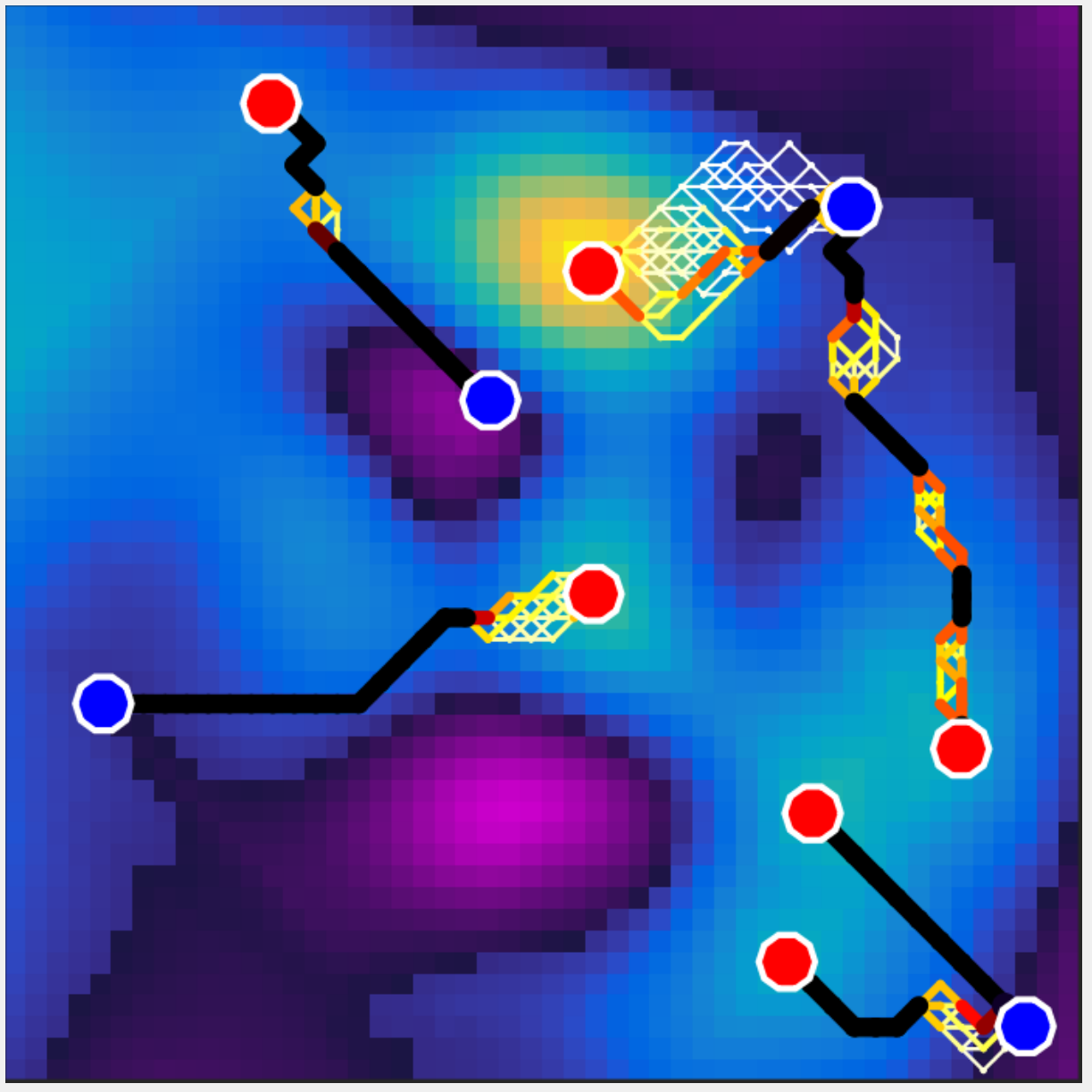}\\
            \vspace*{1mm}
            \includegraphics[scale=0.9]{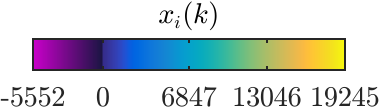}
            \caption{Map of 
            $\bar x_i(k)$, constr. system, $C_{max}{=}25$, stochastic choices.}
            \label{fig:sim2:map:xcon}
        \end{subfigure}
        \caption{(\textit{Example 2}). The map of the network. Red circles: source nodes; blue circles: sink nodes, colored pixels: nodes of the network; color of each pixel: the corresponding value of $h_i$, $\bar{\bar x}_i$ or $\bar x_i(k)$, respectively, with gray representing 0.
        In Fig. \ref{fig:sim2:map:val}, there are hills in the yellow areas, and the altitude decreases as the color becomes bluer: arc costs $\gamma_{ij}$ are negative for the arcs directed downhill, and non-negative otherwise. 
        In Figs. \ref{fig:sim2:map:barx}-\ref{fig:sim2:map:xcon}, there is a higher accumulation of tokens in the yellow areas; the number of tokens is below the zero reference level (negative states) in the violet areas. 
        In Fig. \ref{fig:sim2:map:barx}, by the definition of  $\bar{\bar x}_i$, colors also represent the length of the shortest path from each node to the closer sink.  
        In Figs. \ref{fig:sim2:map:xuncdet}-\ref{fig:sim2:map:xcon},  the lines represent the paths traveled by  $10000$ injected tokens; from white thin lines, when few tokens traversed each arc, to thicker black lines, when almost all the tokens traversed them.}\label{fig:sim2:map}  
        \end{figure}

        We now study a large grid network described by the map of size $50\times50$ pixels represented in Fig. \ref{fig:sim2:map:val}, where each pixel is associated with a node $i$
        and an integer  ``altitude''~$h_i$.
        The network includes all the arcs from each node to the existing nodes in its 8-neighborhood.
        To assign a cost to each arc $(i,j)$, 
        given $dh=h_j-h_i$, we compute: 
        $$
        \gamma_{ij}=
        \begin{cases}
        \text{ceil}(m^-(dh-h_0)) &\textnormal{if~~}dh\leq h_0,\\
        \text{ceil}(m^+(dh-h_0)) &\textnormal{otherwise,}
        \end{cases}
        $$
        where $h_0\in\mathbb{Z}, h_0<0,$ is the difference $dh$ associated with a zero cost, and
        $m^-, m^+\in\mathbb{R}$, $0<m^-<m^+$, ensure positive circuits.
        Negative gradients of $h_i$ such that $dh<h_0$  corresponds to negative cost arcs. We use: $h_0=-30,m^-=0.4$ and $m^+=0.9$. Also, we set $\sigma_{ij}=1$ for each arc.
        The resulting network has $2500$ nodes and $19404$ arcs, with $\gamma_{ij}\in[-97, 273]$. We set six sources $s_1, \dots, s_6$ and four sinks $d_1, \dots, d_4$ (red and blue dots in Fig. \ref{fig:sim2:map:val}, respectively). Fig. \ref{fig:sim2:map:barx} shows the maximal rest state $\bar{\bar x}$ for the unconstrained system, i.e., the minimum distance from each node to the closer sink.
        
        We compare the deterministic and stochastic choice models. The unconstrained system reaches the similar global rest states $\bar x(k)$ shown in Figs.  \ref{fig:sim2:map:xuncdet} and   \ref{fig:sim2:map:xunc}, respectively,  around $k=10^6$. 
        The paths traveled by  $10000$  tokens injected next are also shown, which are the shortest ones to the closest sink. 
        Indeed, in both cases, $\bar x(k)\leq \bar{\bar x}$  and $\bar x_{s_1}=383$, $\bar x_{s_2}=1058$, $\bar x_{s_3}=529$, $\bar x_{s_4}=597$, $\bar x_{s_5}=635$, $\bar x_{s_6}=608$, 
        which are the optimal path lengths (see Table \ref{tab:opt_len}).
        When paths are constrained by $C_{max}=25$, the system reaches the similar global rest states in Figs. \ref{fig:sim2:map:xcondet} and \ref{fig:sim2:map:xcon}, respectively,  around $k=1.35\cdot10^7$ (for simplicity, $\bar x_i(k)=\sum_{c=0}^{C_{max}}\bar x^c_i(k)$ is shown) and then tokens follows the depicted shortest feasible paths.
        The unconstrained shortest paths from $s_1$, $s_5$, $s_6$ are not feasible, 
        so new feasible shortest ones are formed, reaching different sinks.
        Those from $s_2, s_3, s_4$ are feasible and do not change.
        Indeed, 
        $\bar x_{s_1}^{0}=496$,
        $\bar x_{s_2}^{0}=1058$, 
        $\bar x_{s_3}^{0}=529$, 
        $\bar x_{s_4}^{0}=597$, 
        $\bar x_{s_5}^{0}=660$, 
        $\bar x_{s_6}^{0}=640$,
        which are the optimal feasible path lengths.        
       Multiple equivalent outgoing shortest paths exist and are discovered and traversed with a stochastic model, while only one of them is chosen with deterministic choices (the others are still potentially available by changing the choice model).
       This is the main difference between the two choice models,  while the performance are similar, with the stochastic one being only slightly worse. 
        
        \begin{table}[t]
        \addtolength{\tabcolsep}{-1.5pt}
        \centering
        \begin{threeparttable}
        \caption{(\textit{Example 2}). Minimum length $L_{min}$ and minimum secondary cost $C_{min}$ of the paths between each source $s_i$ and sink $d_i$.}
        \label{tab:opt_len}
        \footnotesize
        \begin{tabular}{ccccccccc}
        \toprule
                       & \multicolumn{2}{c}{\textbf{$d_1$}}            & \multicolumn{2}{c}{\textbf{$d_2$}} & \multicolumn{2}{c}{\textbf{$d_3$}}           & \multicolumn{2}{c}{\textbf{$d_4$}} \\
        \cmidrule(lr){2-3}\cmidrule(lr){4-5}\cmidrule(lr){6-7}\cmidrule(lr){8-9}
        & \textit{$L_{min}$}             & \textit{$C_{min}$}           & \textit{$L_{min}$}       & \textit{$C_{min}$}      & \textit{$L_{min}$}            & \textit{$C_{min}$}           & \textit{$L_{min}$}       & \textit{$C_{min}$}      \\
        \midrule
        \textbf{$s_1$} & 608                    & 27                   & \textit{496}     & \textit{14}     & \textbf{383}          & \textbf{28}          & 1527             & 43              \\
        \textbf{$s_2$} & \textit{\textbf{1058}} & \textit{\textbf{12}} & 1141             & 6               & 1367                  & 23                   & 2133             & 35              \\
        \textbf{$s_3$} & 637                    & 18                   & 637              & 9              & \textit{\textbf{529}} & \textit{\textbf{23}} & 1145             & 20              \\
        \textbf{$s_4$} & \textit{\textbf{597}}  & \textit{\textbf{25}} & 817              & 22              & 776                   & 40                   & 720              & 13              \\
        \textbf{$s_5$} & 698                    & 28                   & 750              & 19              & \textbf{635}          & \textbf{33}          & \textit{660}     & \textit{10}     \\
        \textbf{$s_6$} & 844                    & 35                   & 896              & 26              & \textbf{608}          & \textbf{32}          & \textit{640}     & \textit{11}    \\
        \bottomrule
        \end{tabular}
        \begin{tablenotes}[flushleft, para]\scriptsize
        \note{shortest paths have secondary cost  $C_{min}$ (it is $31$ for $s_2\rightarrow d_3$). For each source, values in bold refer to the shortest (unconstrained) outgoing path, and those in italic to the shortest feasible one when $C_{max}{=}25$.}
        \end{tablenotes}
        \end{threeparttable}
        \end{table}
        
\begin{figure}[t]
        \centering
        \begin{subfigure}[t]{0.24\linewidth}
            \includegraphics[width=1\linewidth]{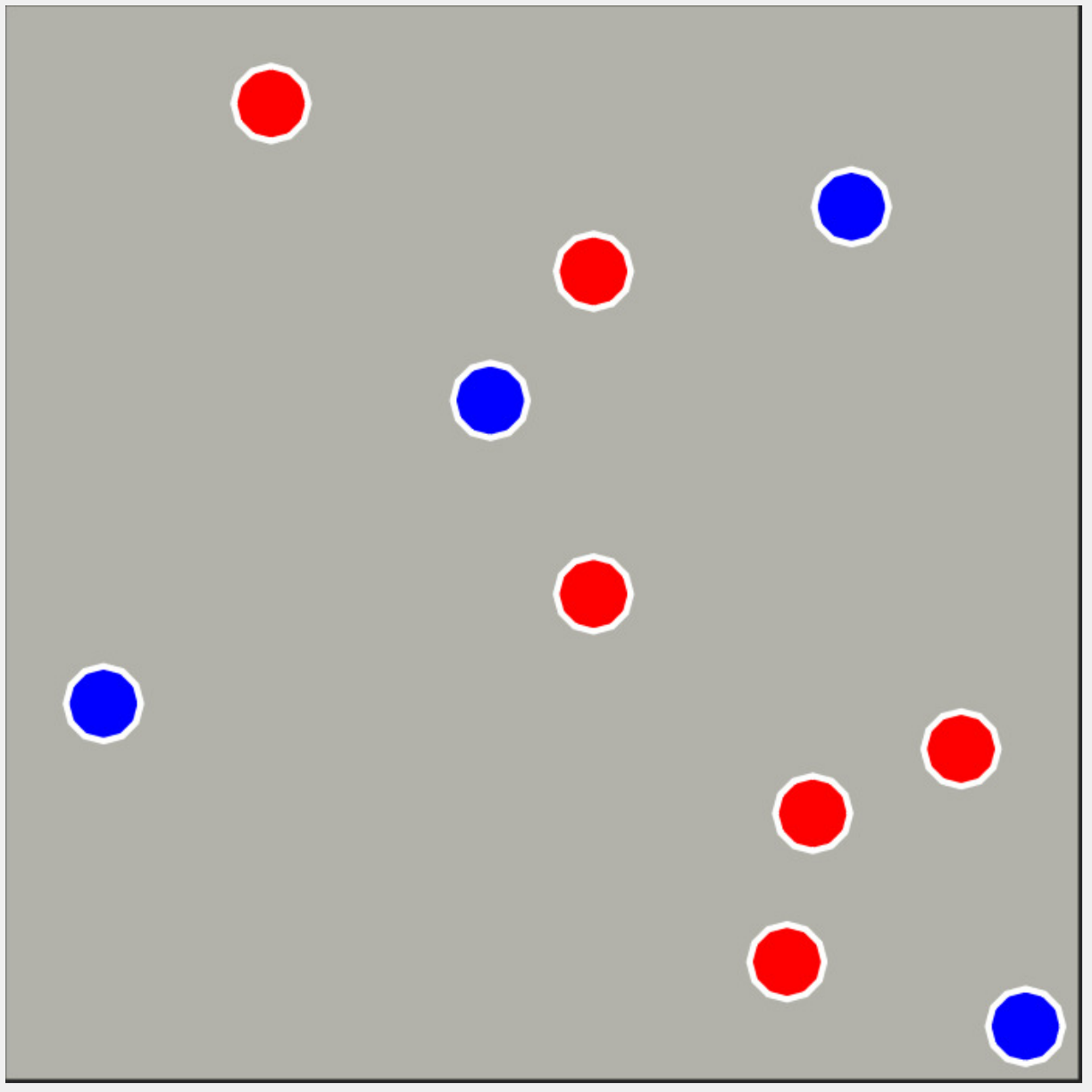}
            \caption{$k=0$.}
            \label{fig:sim2:tran:0}
        \end{subfigure}
        \hfill
        \begin{subfigure}[t]{0.24\linewidth}
            \includegraphics[width=1\linewidth]{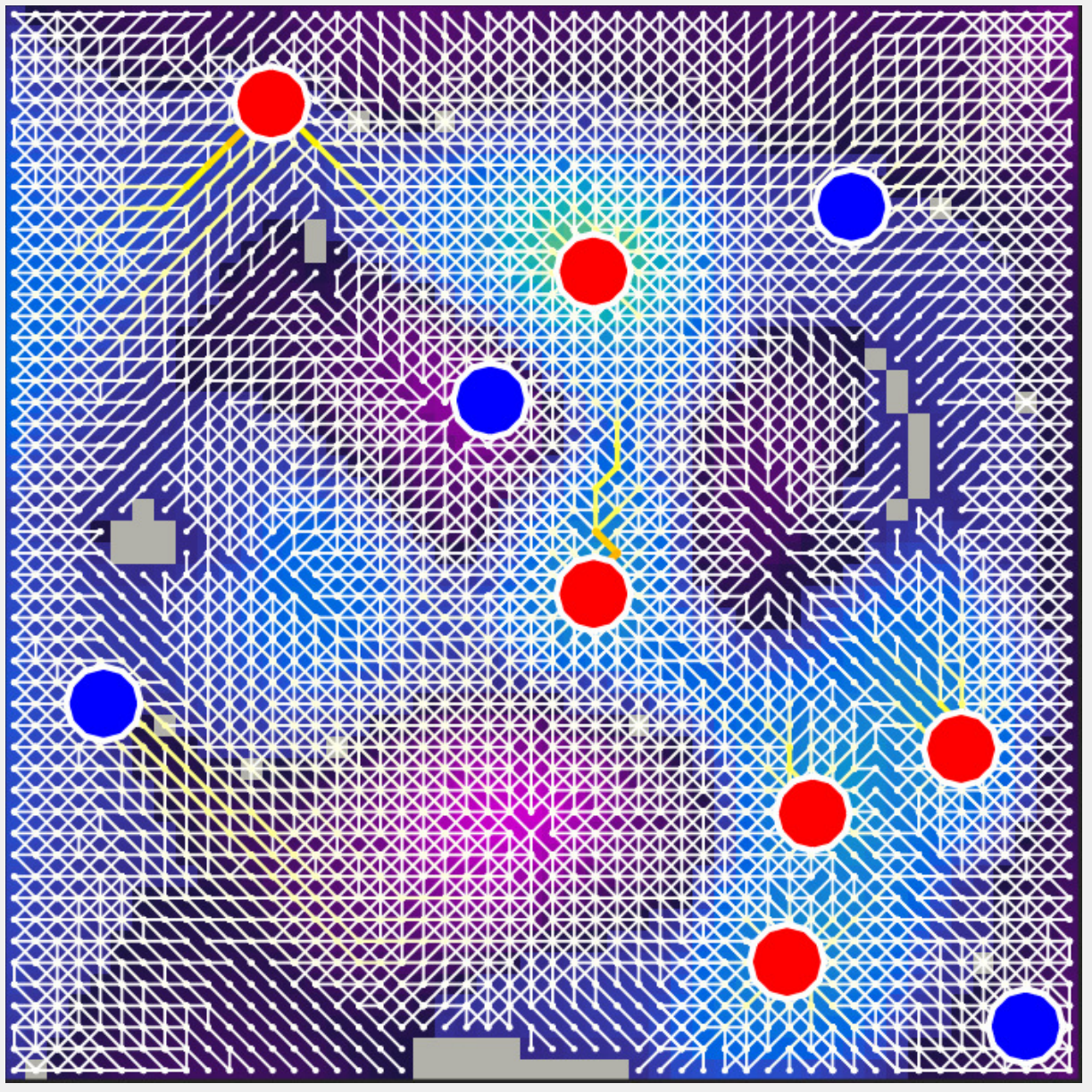}
            \caption{$k=1\cdot10^5$.}
            \label{fig:sim2:tran:1}
        \end{subfigure}
        \hfill
        \begin{subfigure}[t]{0.24\linewidth}
            \includegraphics[width=1\linewidth]{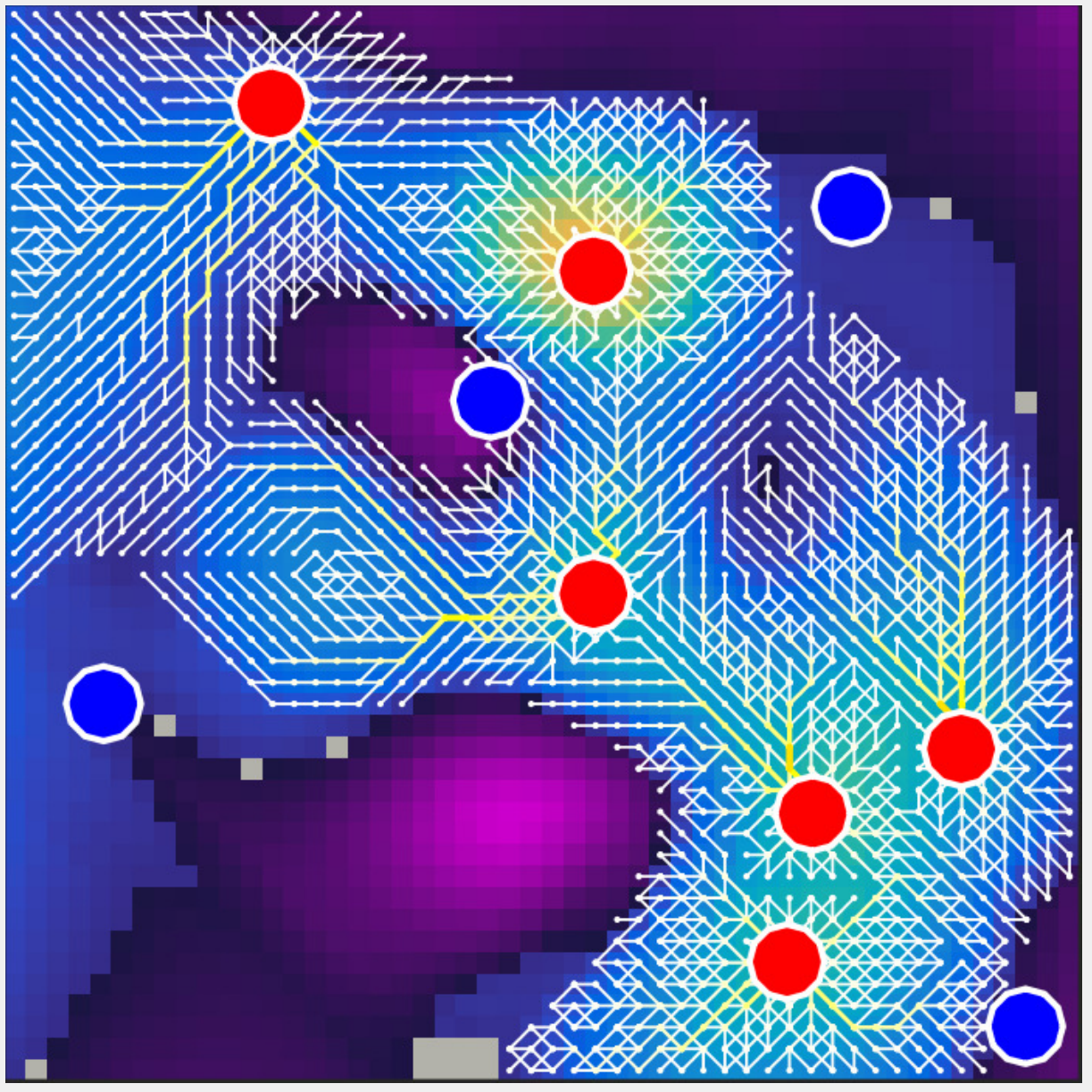}
            \caption{$k=3\cdot10^5$.}
            \label{fig:sim2:tran:2}
        \end{subfigure}
        \hfill
        \begin{subfigure}[t]{0.24\linewidth}
            \includegraphics[width=1\linewidth]{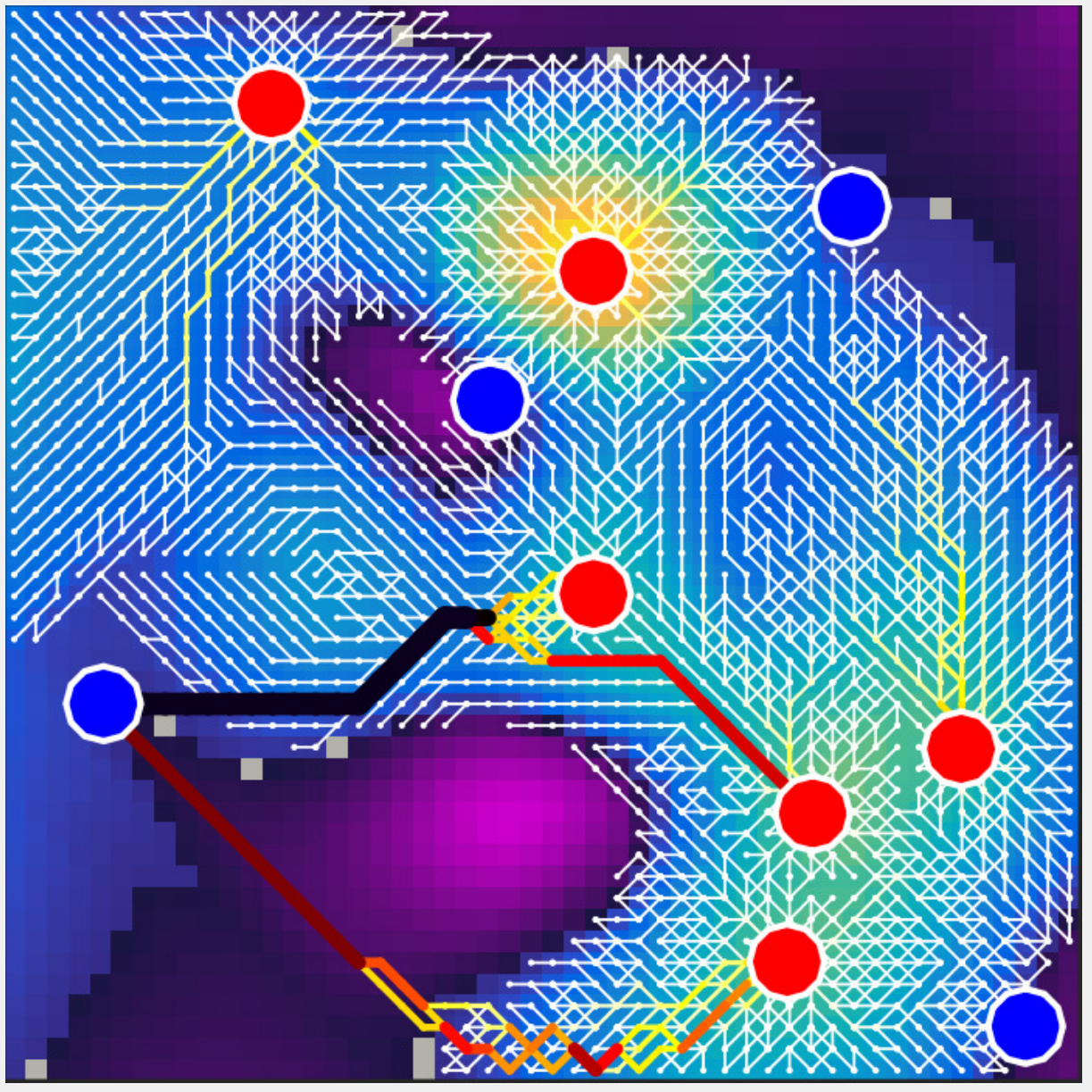}
            \caption{$k=5\cdot10^5$.}
            \label{fig:sim2:tran:3}
        \end{subfigure}\\
        \vspace*{2mm}
        \begin{subfigure}[t]{0.24\linewidth}
            \includegraphics[width=1\linewidth]{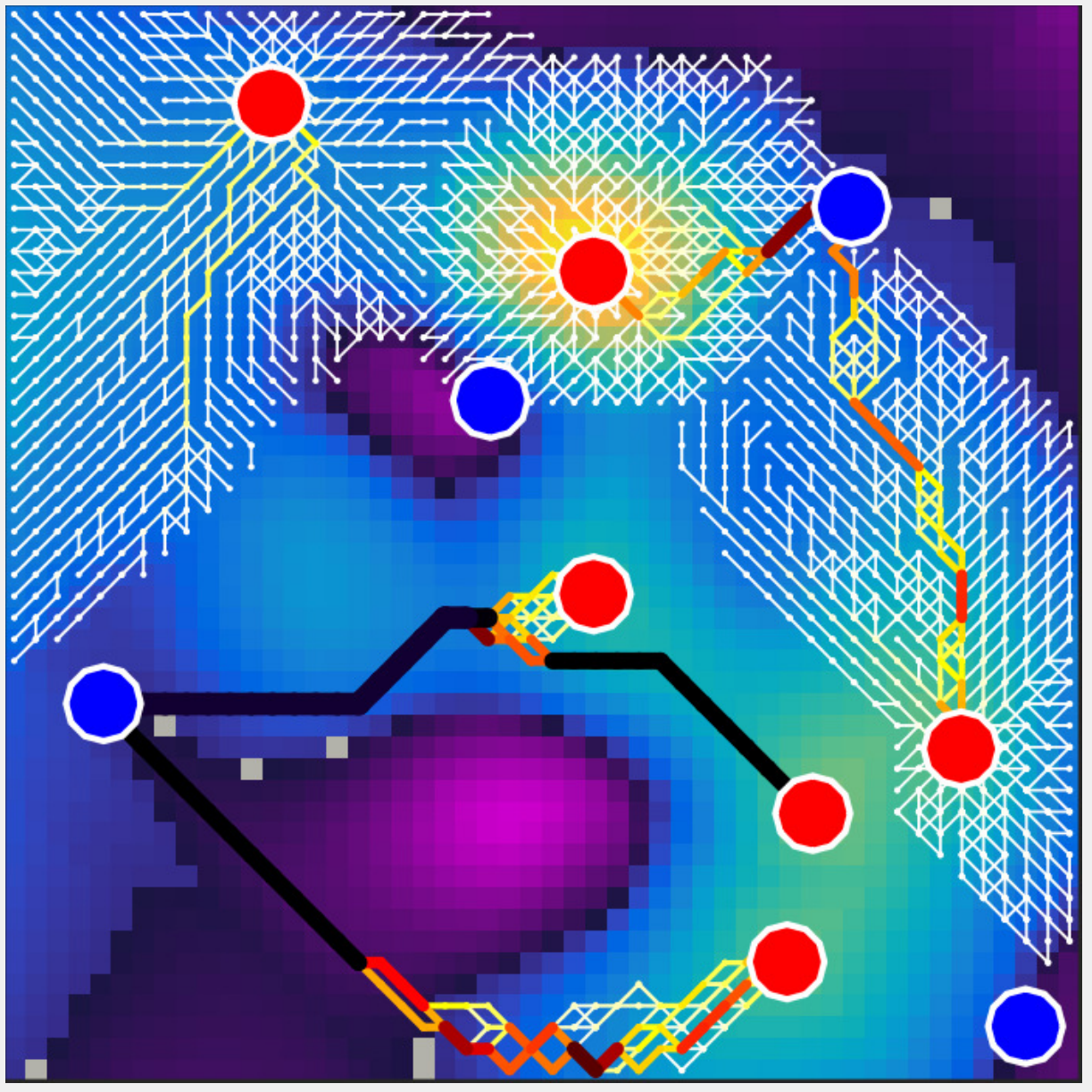}
            \caption{$k=6\cdot10^5$.}
            \label{fig:sim2:tran:4}
        \end{subfigure}
        \hfill    
        \begin{subfigure}[t]{0.24\linewidth}
            \includegraphics[width=1\linewidth]{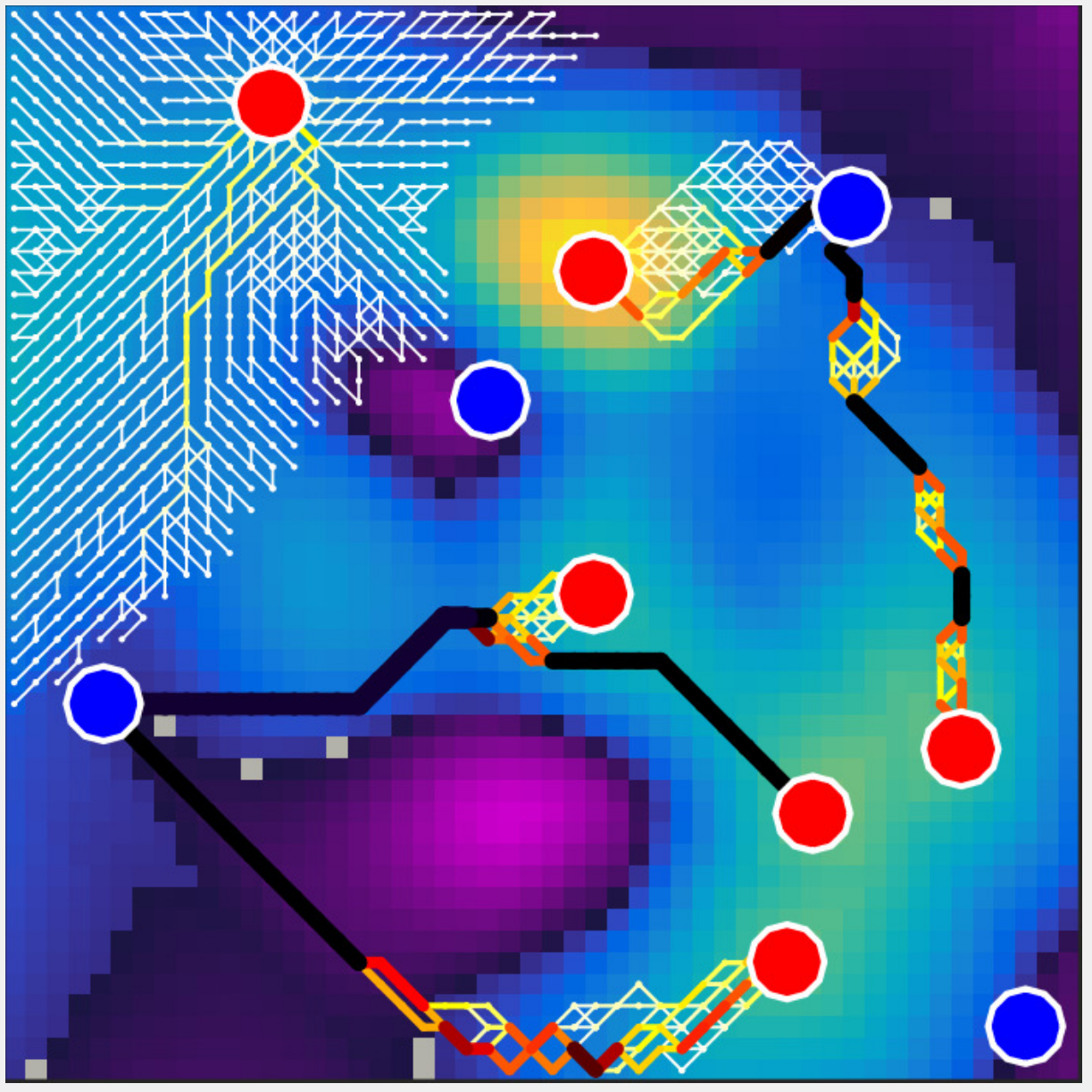}
            \caption{$k=8\cdot10^5$.}
            \label{fig:sim2:tran:5}
        \end{subfigure}
        \hfill
        \begin{subfigure}[t]{0.24\linewidth}
            \includegraphics[width=1\linewidth]{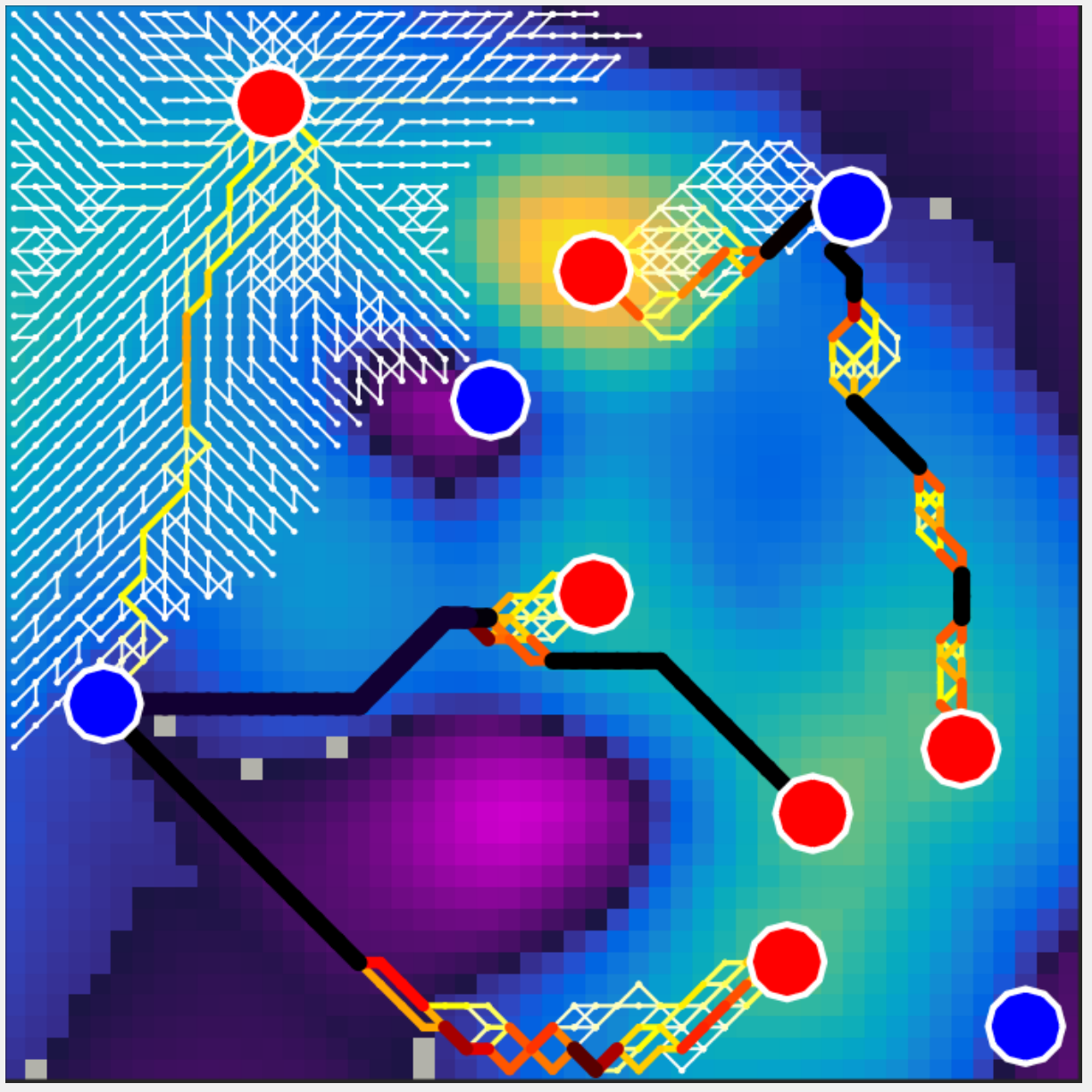}
            \caption{$k=10\cdot10^5$.}
            \label{fig:sim2:tran:6}
        \end{subfigure}
        \hfill
        \begin{subfigure}[t]{0.24\linewidth}
            \includegraphics[width=1\linewidth]{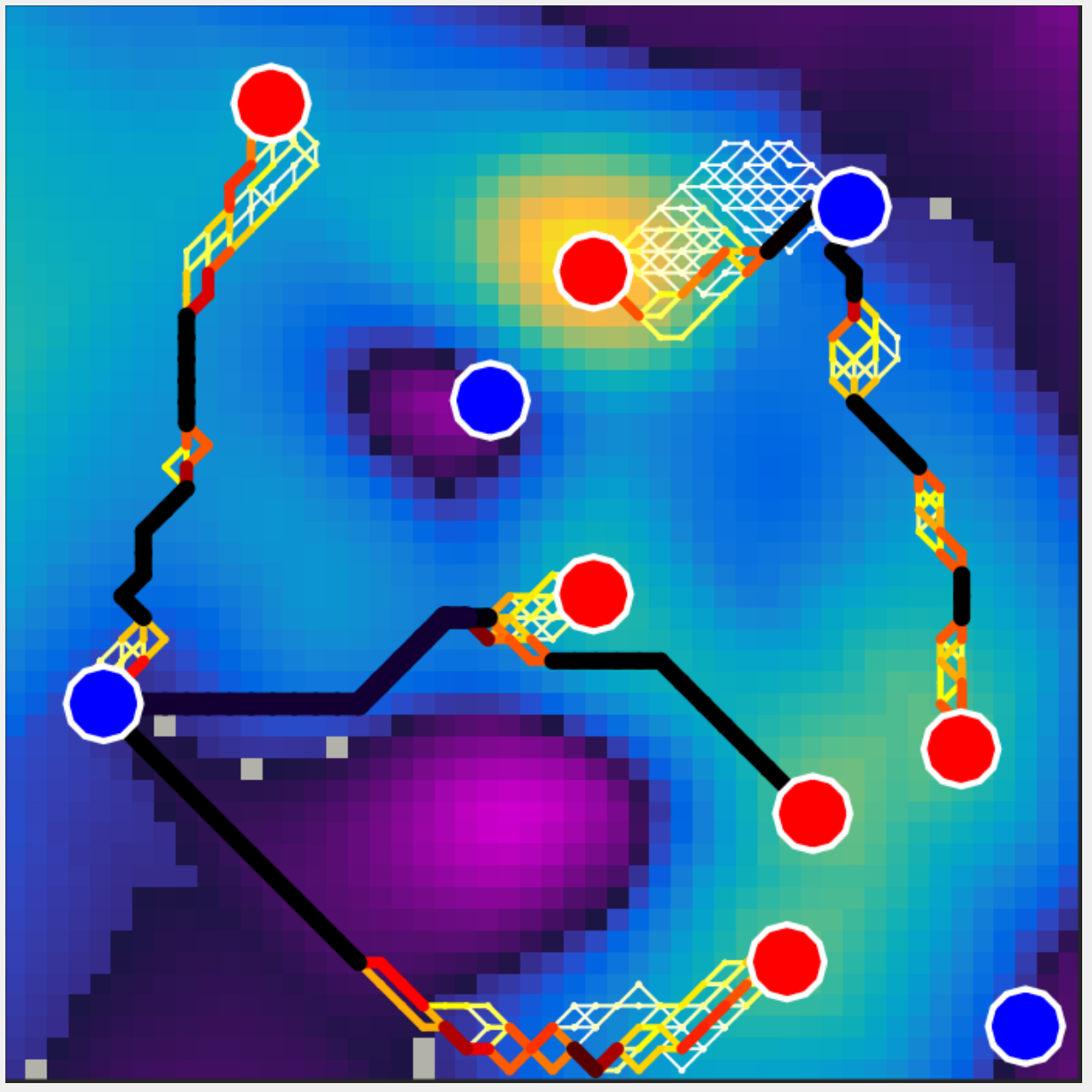}
            \caption{$k=12\cdot10^5$.}
            \label{fig:sim2:tran:7}
        \end{subfigure}
        \hfill\\
        \caption{(\textit{Example 2}). Evolution of the nodes' state map of the unconstrained system, to get the final state of Fig. \ref{fig:sim2:map:xunc}. The paths traveled by the tokens in the previous $10^5$ time units are represented.}
        \label{fig:sim2:tran}
        \end{figure}

        The transient is shown in Fig. \ref{fig:sim2:tran} for the stochastic unconstrained policy.
        Initially, the negative cost arcs make the state not admissible. 
        The fast dynamic makes it admissible soon, at time $k=6$: tokens moves simultaneously in different parts of the network along negative cost arcs, i.e., the states of the nodes in the ``hills'' become negative, transferring tokens downhill (such paths are included in Fig. \ref{fig:sim2:tran:1}). Then, only the injected tokens move, spreading in the network and increasing its state by accumulating in its nodes. Their exploration area grows over time until the closer sink 
        is reached; afterwards, they all reach such sink along the shortest paths.
        To get to this point, $552628$ injected tokens have stopped in the nodes.
        The behaviour is similar when paths are constrained, however now the exploration area is restricted; also, as now the number of states is multiplied by $C_{max}+1=26$, 
        the constrained system takes longer to reach a global rest state and more tokens need to accumulate in the nodes: the state becomes admissible at $k=50$ and, in the end, $6134934$ injected tokens have stopped in the nodes, possibly falling asleep.

        After the network stabilizes, we apply the  modifications described by the top rows of Figs. \ref{fig:sim3:mod1}-\ref{fig:sim3:mod6}.  
        Specifically, at $k=18000001$, some nodes  are removed.
        At $k=21000001$, some of these are re-enabled. 
        At $k=24000001$, some nodes stop being sources.
        At $k=27000001$, some nodes become new sources.
        At $k=42000001$, some nodes become new sinks.
        At $k=45000001$, some nodes stop being sinks.
        The system adapts to each change, reaching a new global rest state and updating the paths followed by the injected tokens, which eventually are the outgoing shortest feasible ones for the unconstrained and constrained system, see the middle and bottom rows of Figs. \ref{fig:sim3:mod1}-\ref{fig:sim3:mod6}.  
        
\begin{figure}[t!]
        \centering
        \includegraphics[scale=1]{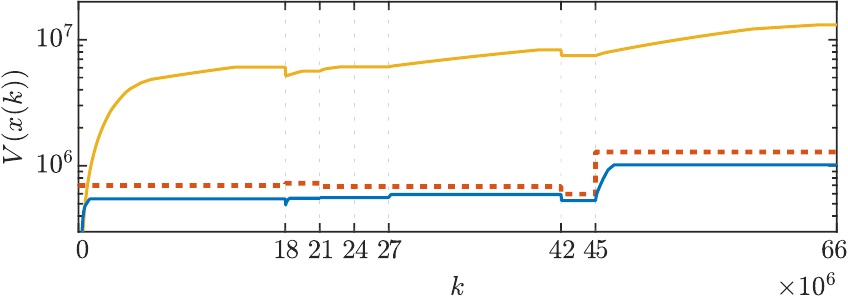}
        \caption{(\textit{Example 2}). Blue: $V(x(k))$ for the unconstrained system; dotted red: $V(\bar{\bar x}(k))$ for the unconstrained system; yellow: $V(x(k))$ for the constrained system with $C_{max}{=}25$. The y-axis is in log-scale.} 
        \label{fig:sim3:time_ev:v} 
        \end{figure}

\begin{figure*}[t]
        \centering
        \begin{subfigure}[t]{0.15\linewidth}
            \includegraphics[width=1\linewidth]{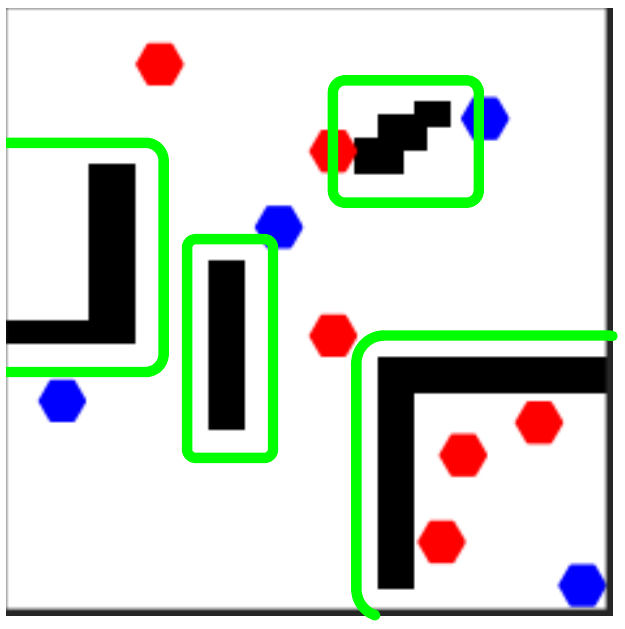}
            \vspace*{1mm}
            \includegraphics[width=1\linewidth]{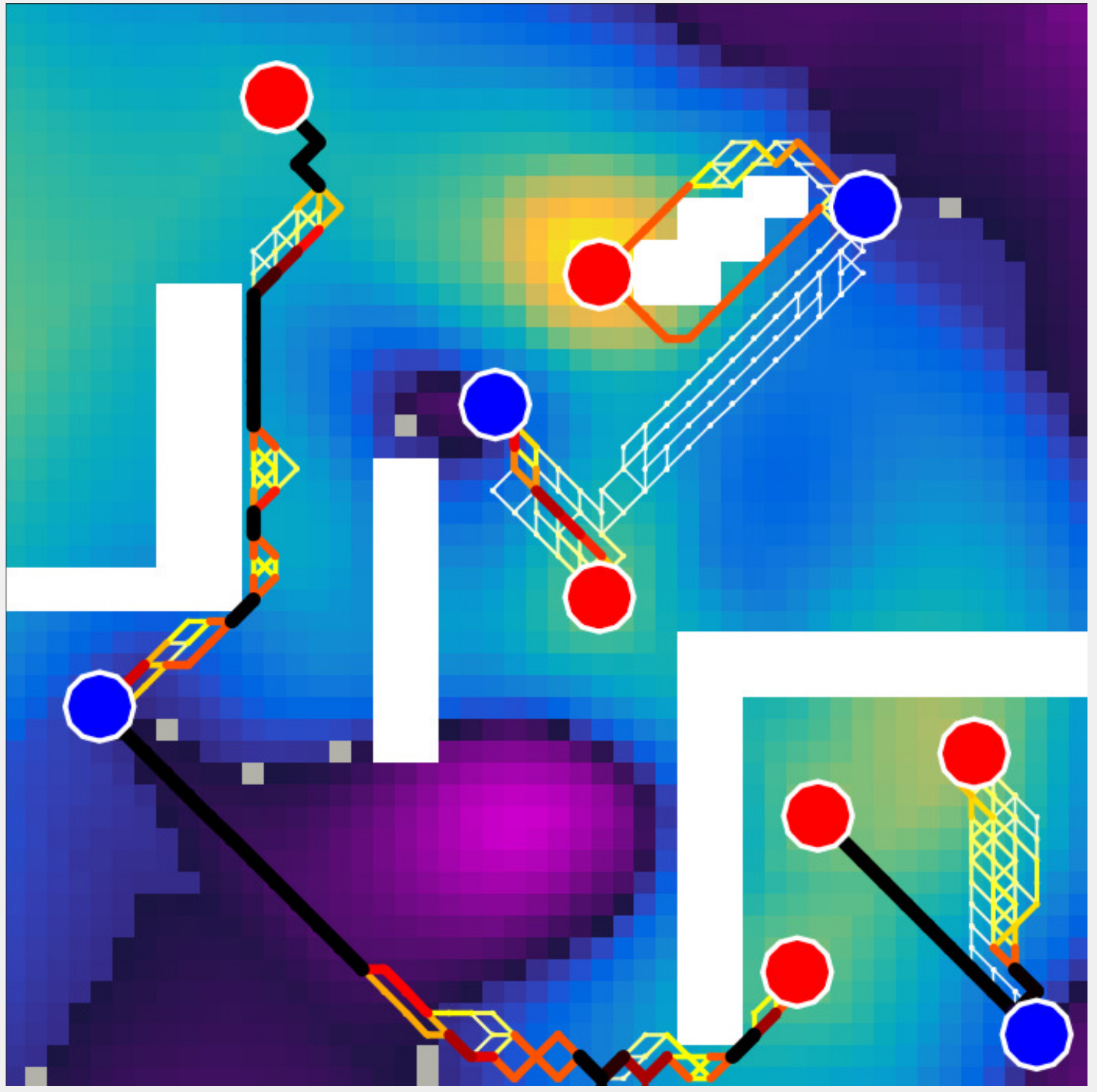}
            \vspace*{1mm}
            \includegraphics[width=1\linewidth]{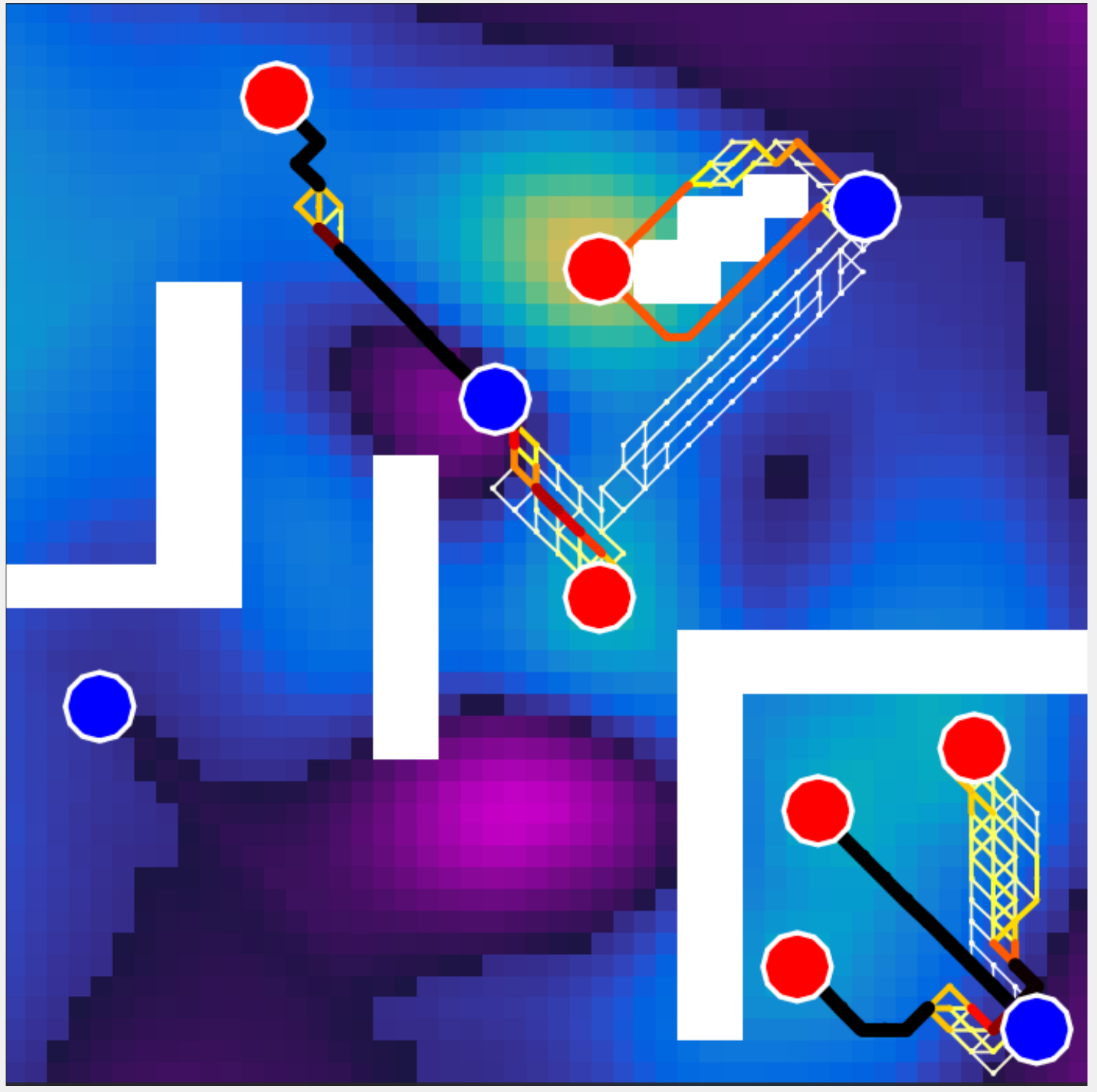}
            \caption{$k_m=18000001$.}
            \label{fig:sim3:mod1}
        \end{subfigure}
        \hfill
        \begin{subfigure}[t]{0.15\linewidth}
            \includegraphics[width=1\linewidth]{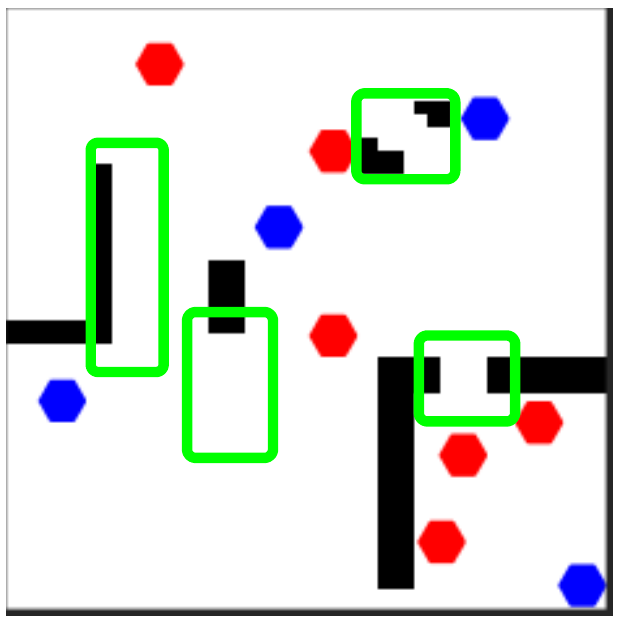}
            \vspace*{1mm}
            \includegraphics[width=1\linewidth]{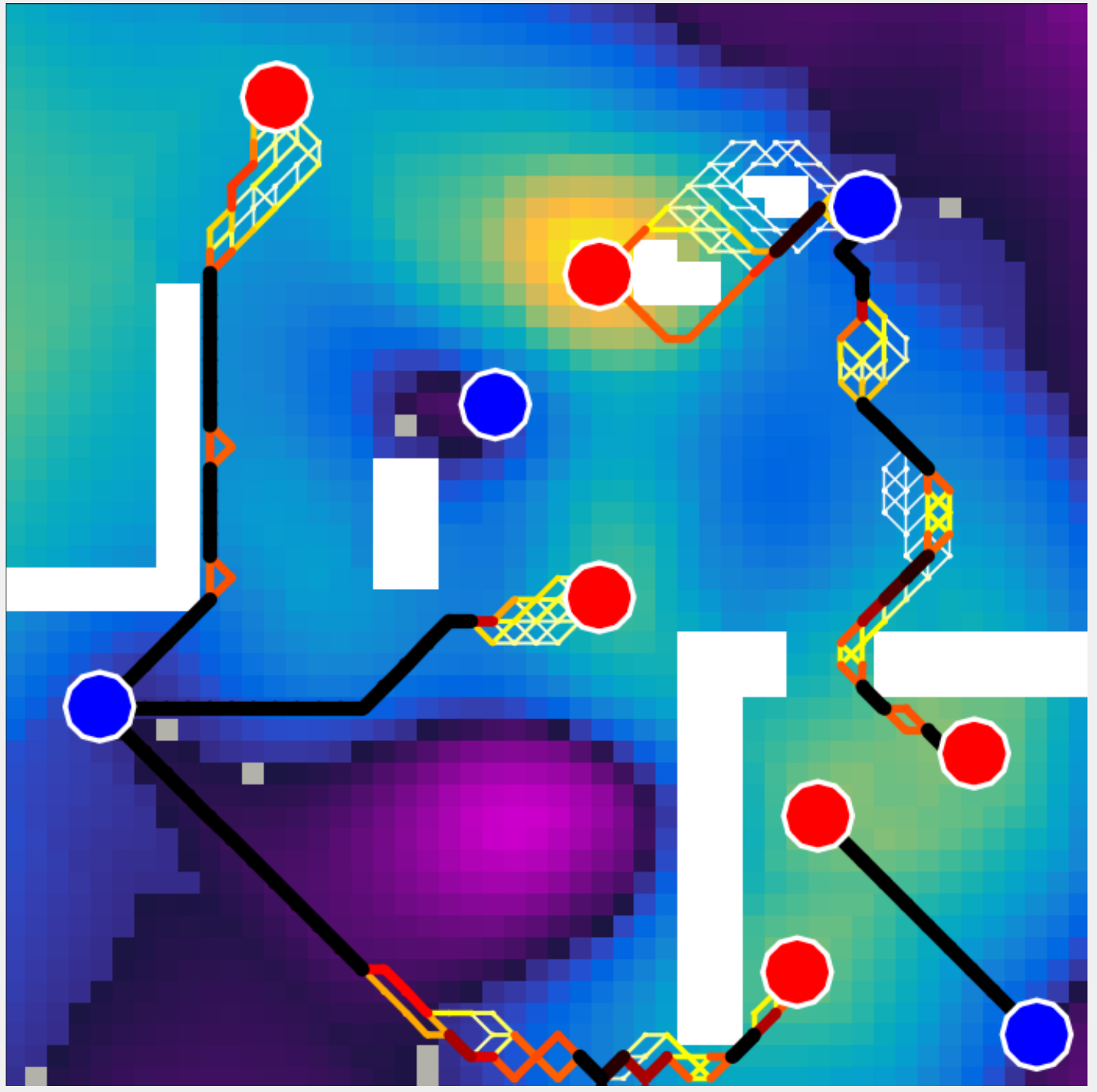}
            \vspace*{1mm}
            \includegraphics[width=1\linewidth]{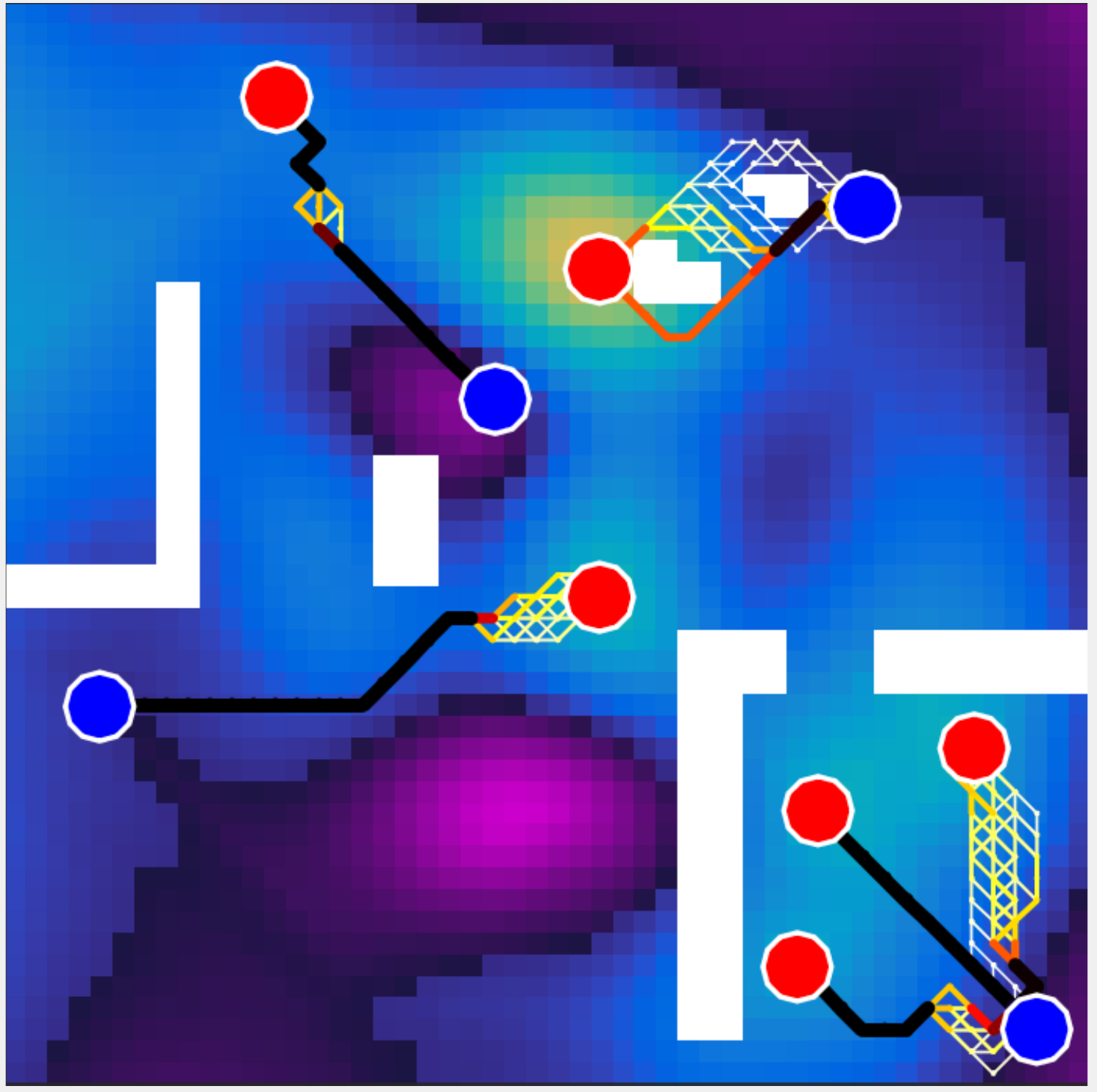}
            \caption{$k_m=21000001$.}
            \label{fig:sim3:mod2}
        \end{subfigure}
        \hfill
        \begin{subfigure}[t]{0.15\linewidth}
            \includegraphics[width=1\linewidth]{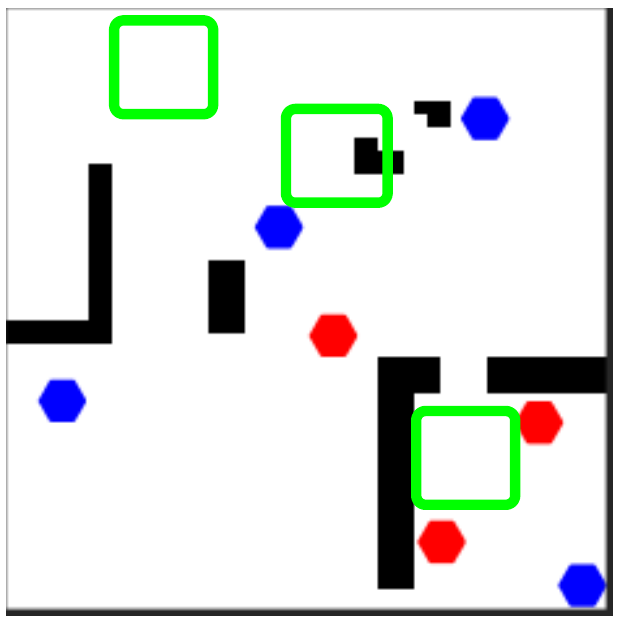}
            \vspace*{1mm}
            \includegraphics[width=1\linewidth]{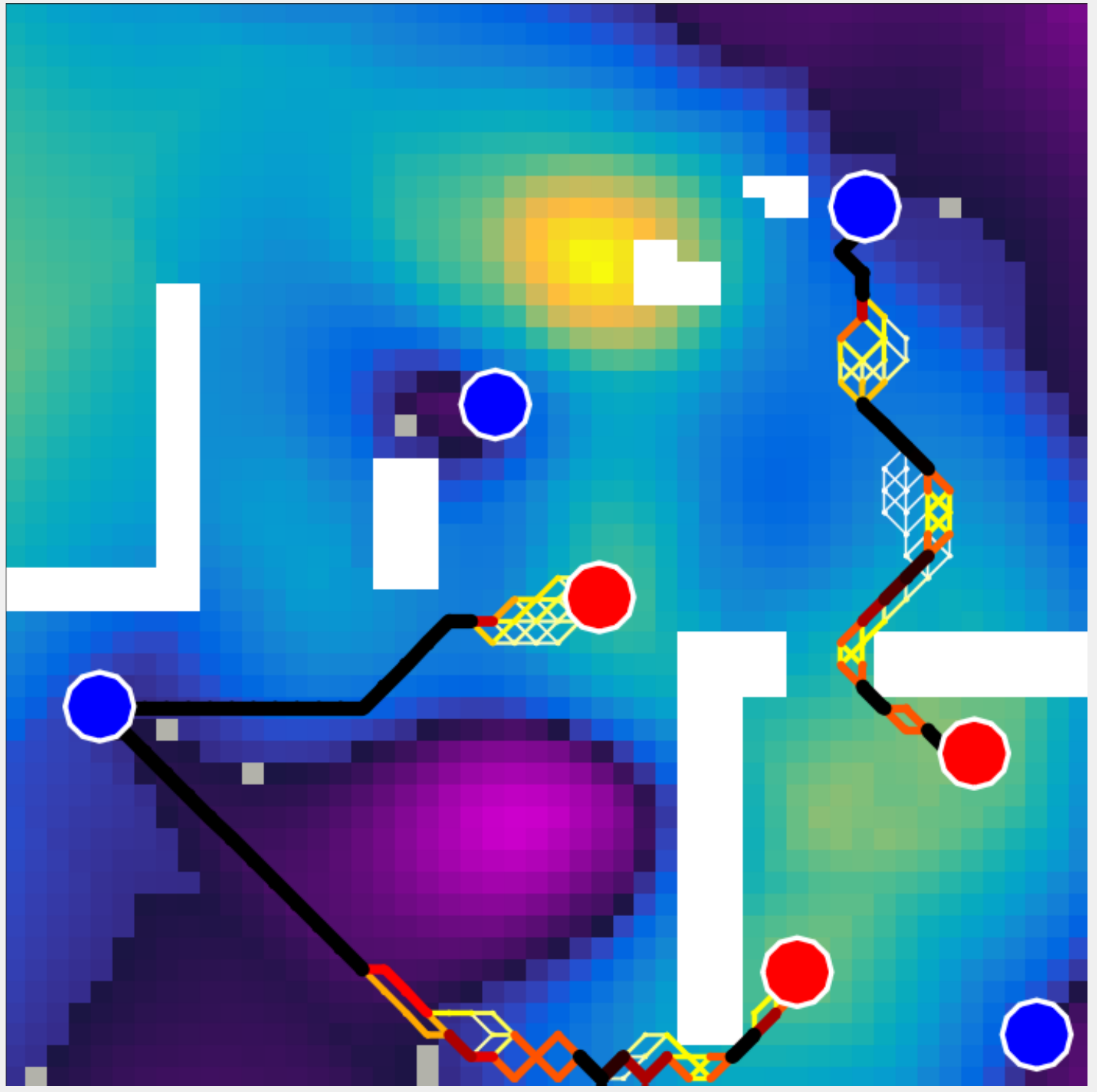}
            \vspace*{1mm}
            \includegraphics[width=1\linewidth]{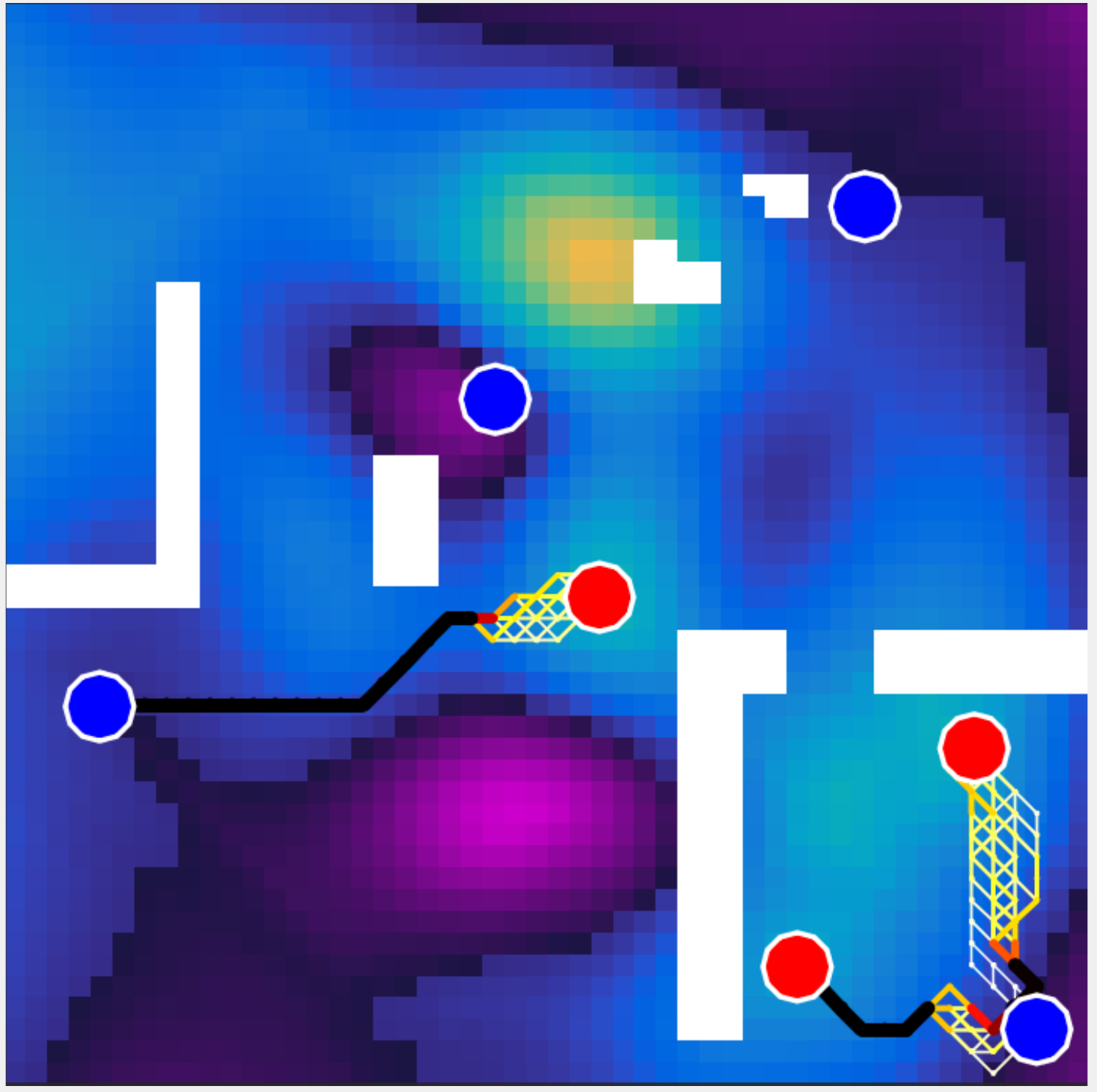}
            \caption{$k_m=24000001$.}
            \label{fig:sim3:mod3}
        \end{subfigure}
        \hfill
        \begin{subfigure}[t]{0.15\linewidth}
            \includegraphics[width=1\linewidth]{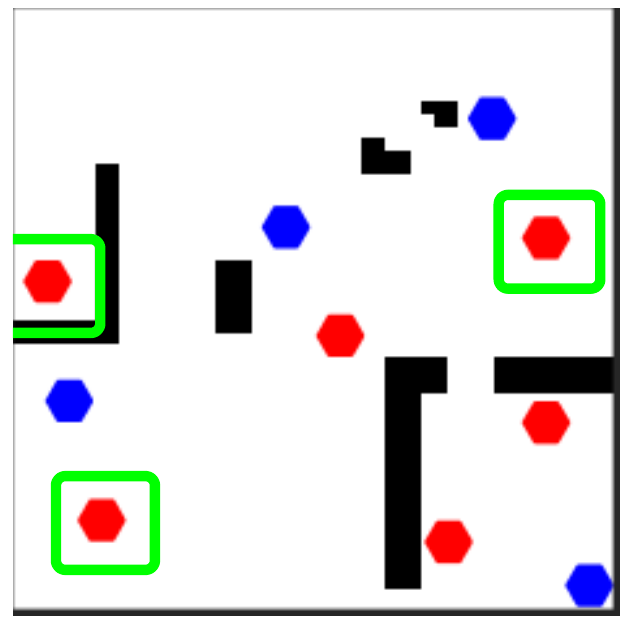}
            \vspace*{1mm}
            \includegraphics[width=1\linewidth]{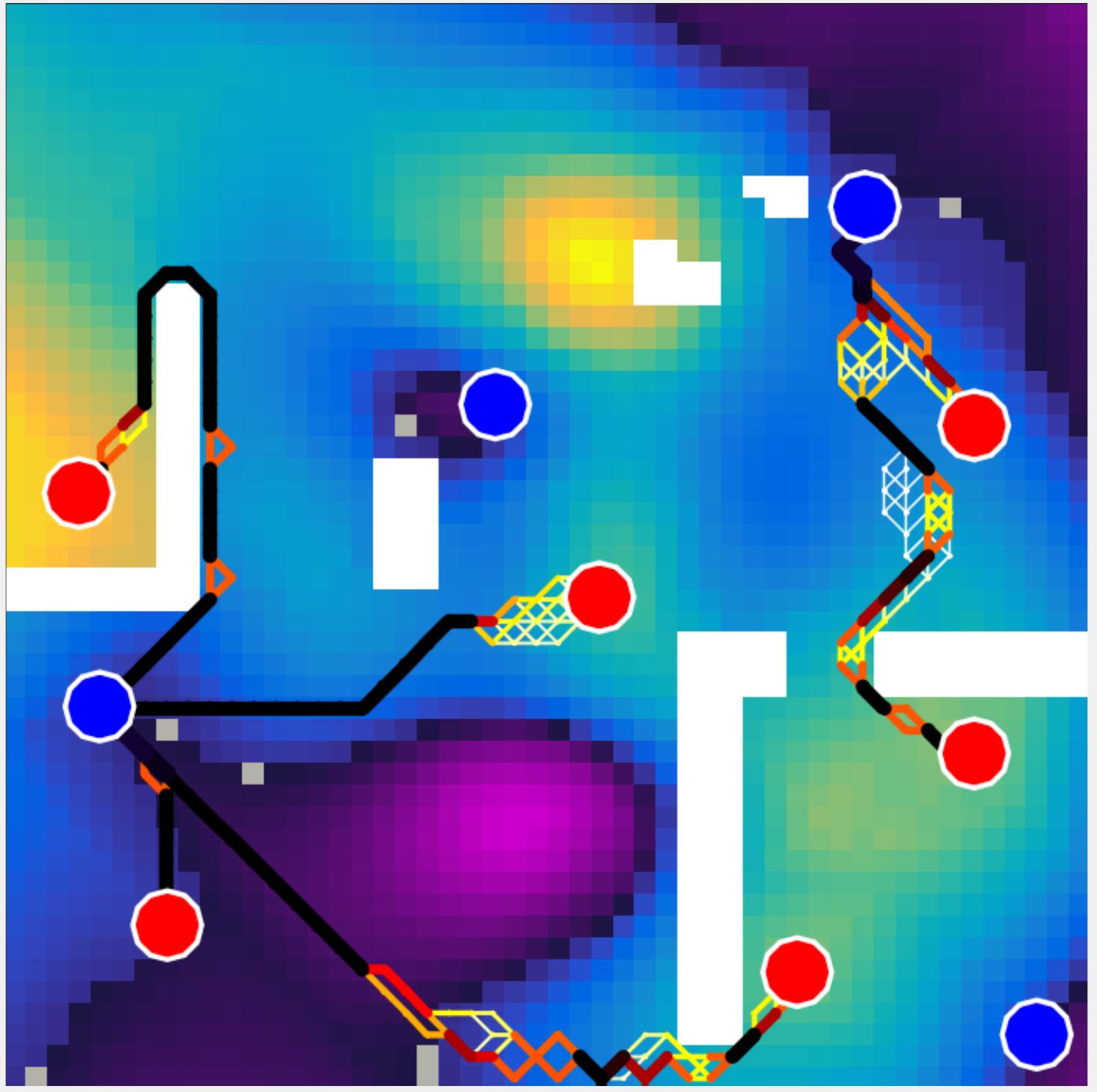}
            \vspace*{1mm}
            \includegraphics[width=1\linewidth]{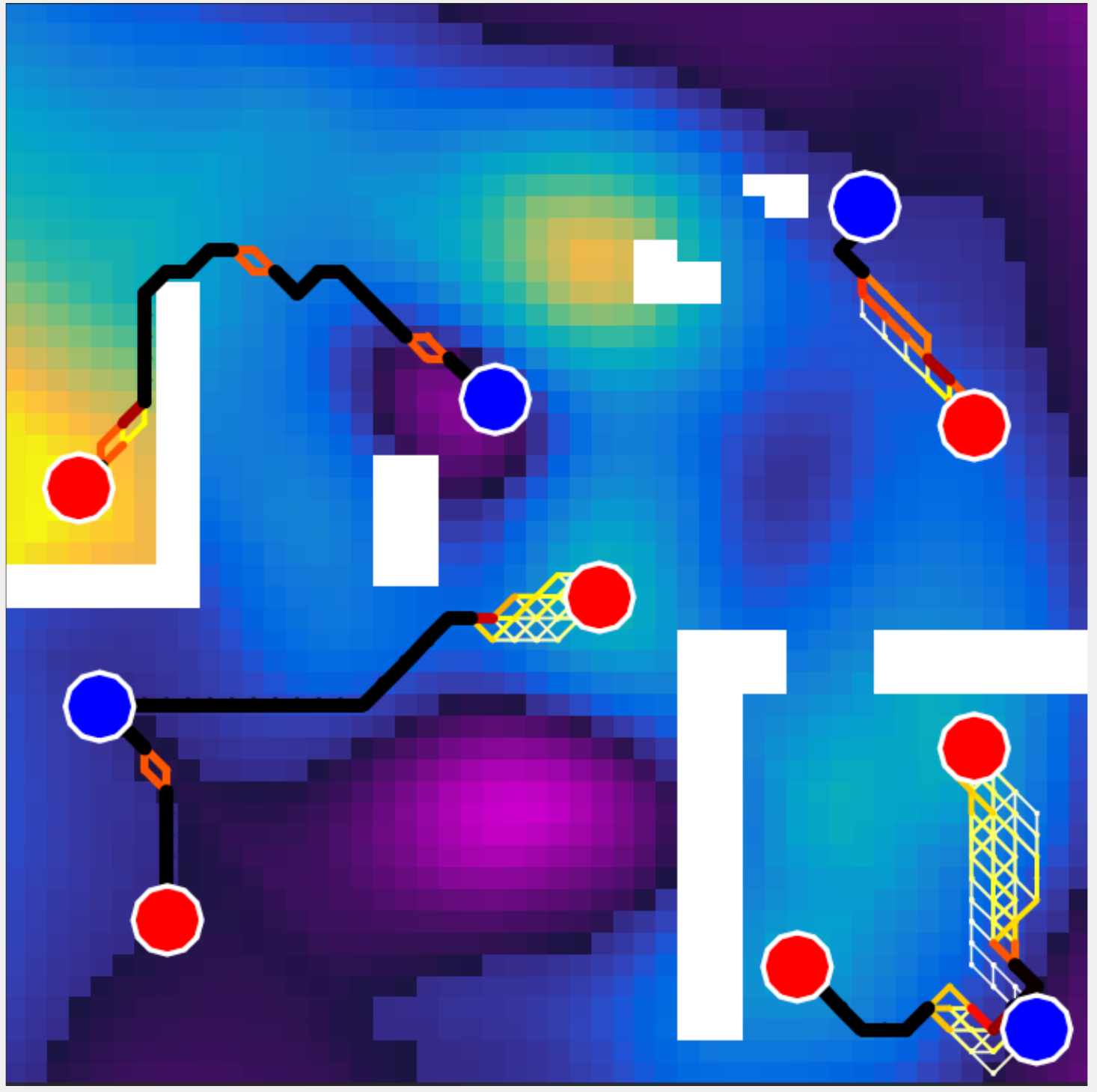}
            \caption{$k_m=27000001$.}
            \label{fig:sim3:mod4}
        \end{subfigure}
        \hfill
        \begin{subfigure}[t]{0.15\linewidth}
            \includegraphics[width=1\linewidth]{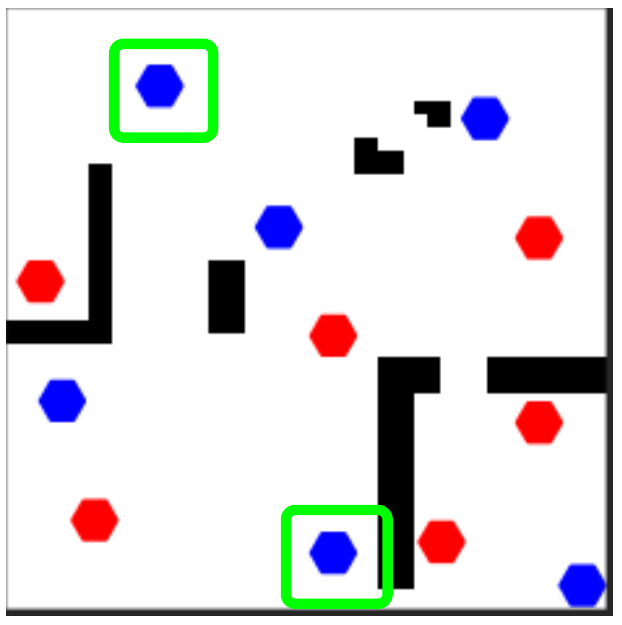}
            \vspace*{1mm}
            \includegraphics[width=1\linewidth]{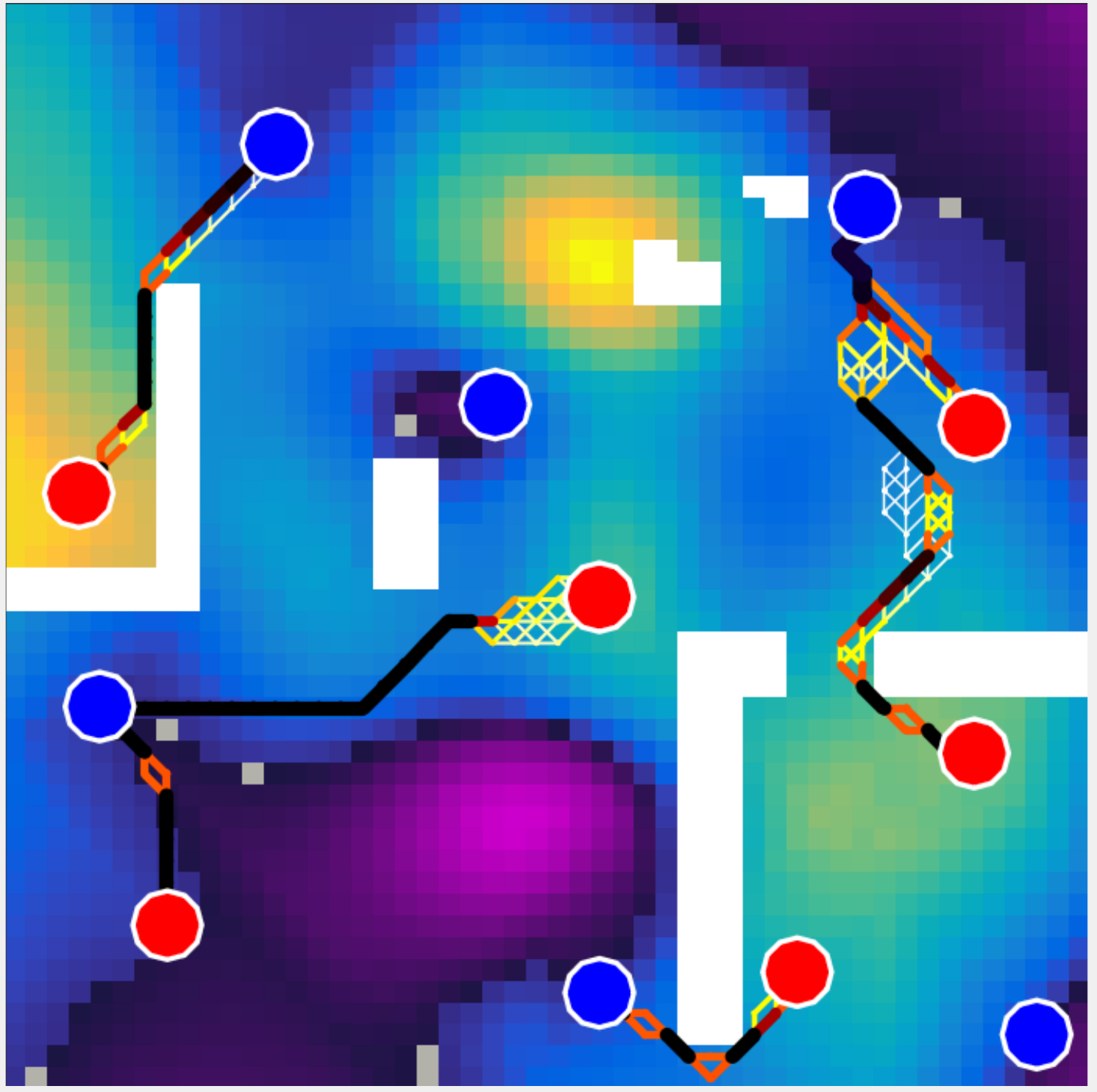}
            \vspace*{1mm}
            \includegraphics[width=1\linewidth]{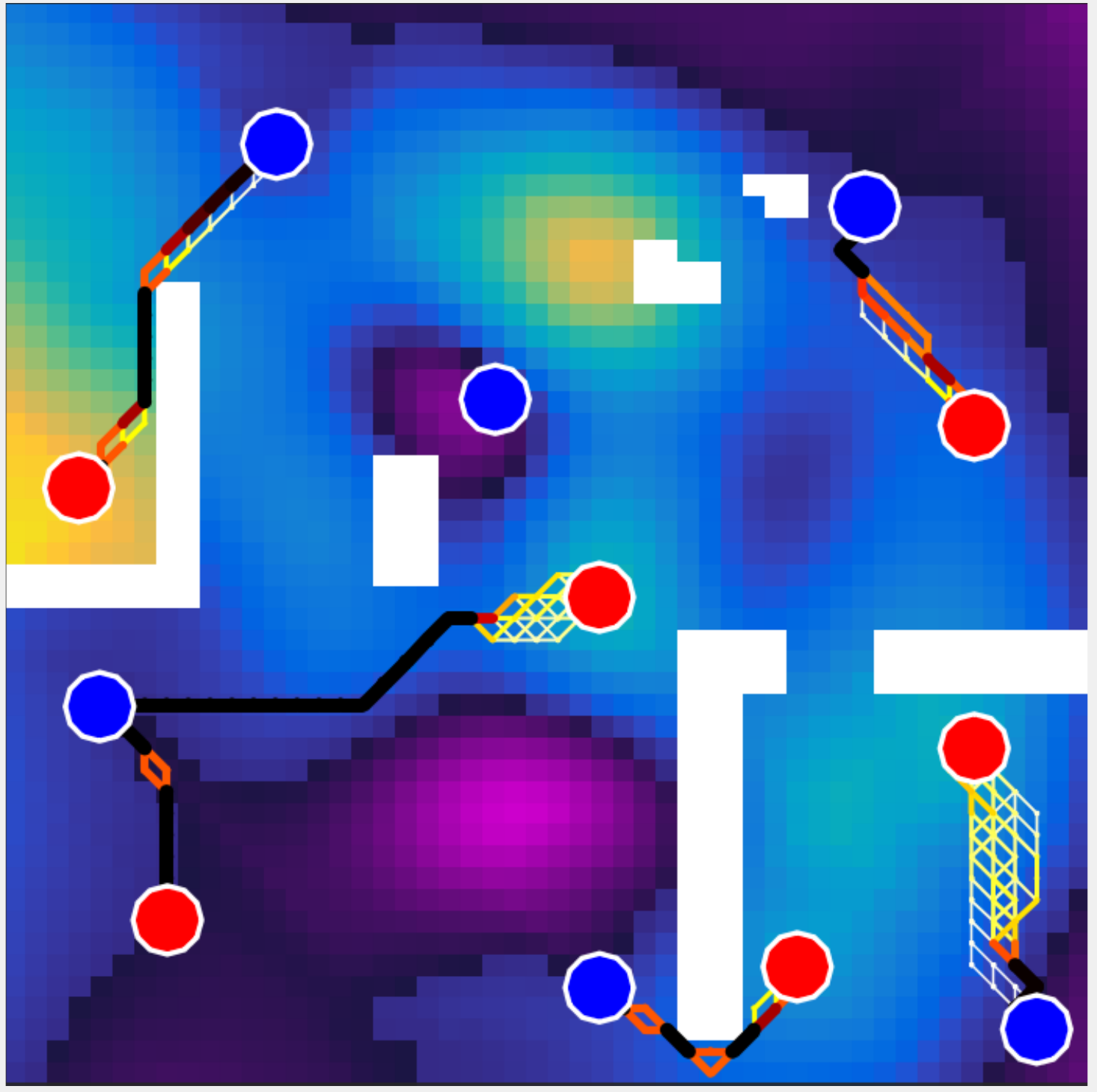}
            \caption{$k_m=42000001$.}
            \label{fig:sim3:mod5}
        \end{subfigure}
        \hfill
        \begin{subfigure}[t]{0.15\linewidth}
            \includegraphics[width=1\linewidth]{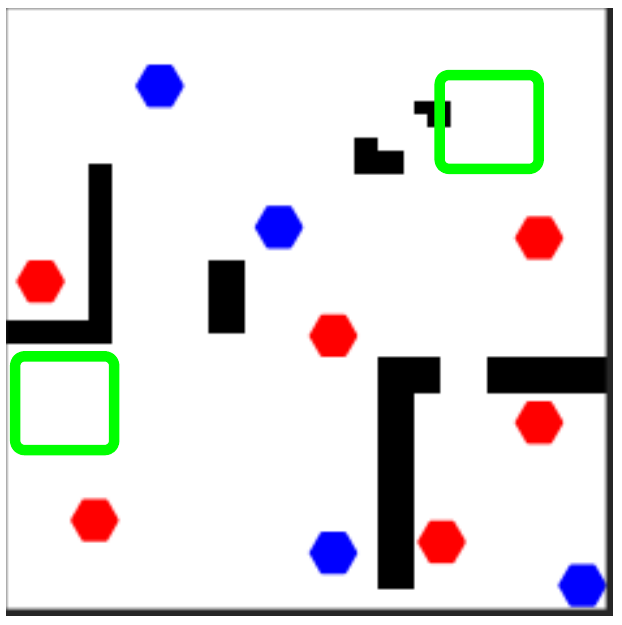}
            \vspace*{1mm}
            \includegraphics[width=1\linewidth]{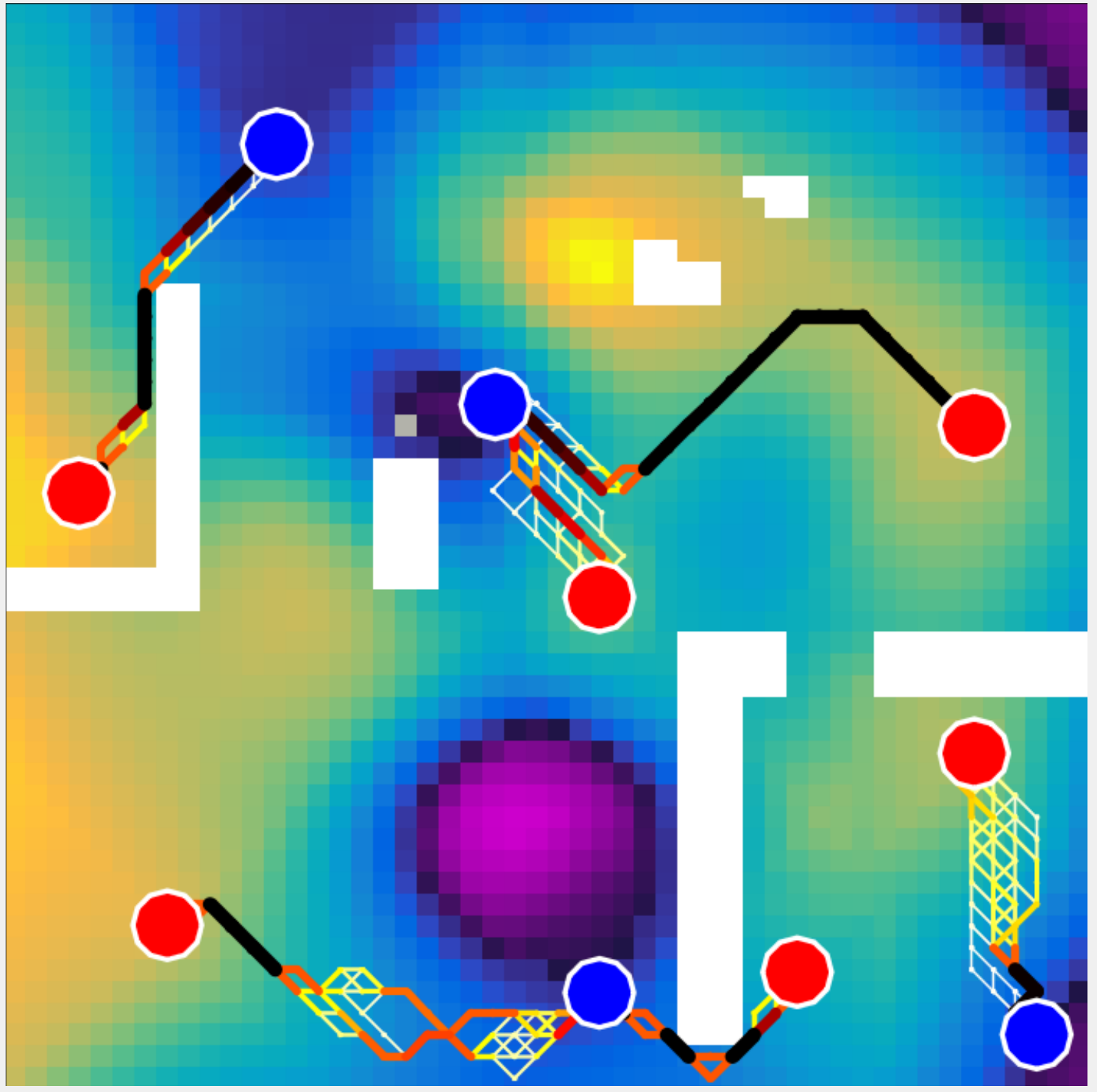}
            \vspace*{1mm}
            \includegraphics[width=1\linewidth]{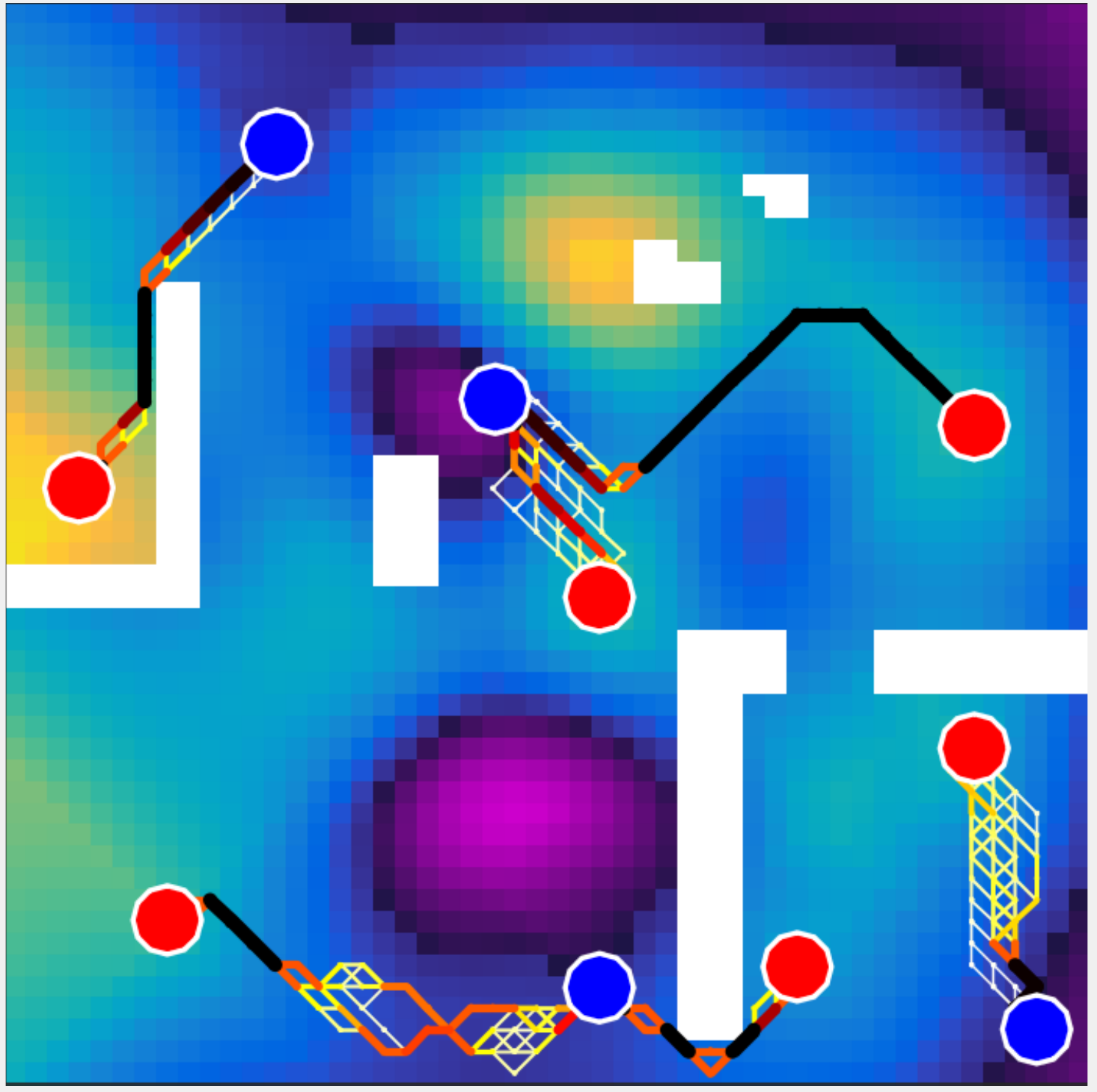}
            \caption{$k_m=45000001$.}
            \label{fig:sim3:mod6}
        \end{subfigure}
        \hfill\\
        \caption{(\textit{Example 2}). \textbf{Top rows:} 
        evolution of the map of Fig. \ref{fig:sim2:map:val}, at different times $k=k_m$. Black pixels: obstacles (disabled nodes); white pixels: enabled nodes; red circles: sources; blue circles: sinks; green shapes: modified sources/sinks/group of nodes. \textbf{Middle and bottom rows:}
        the corresponding global rest state $\bar x_i$ reached after each modification, for the unconstrained (middle row, evolving from Fig. \ref{fig:sim2:map:xunc}) and constrained system with $C_{max}=25$ (bottom row, evolving from Fig. \ref{fig:sim2:map:xcon}). White pixels: obstacles. The paths traveled by $40000$ injected tokens are shown.}
        \label{fig:sim3:mod}
        \end{figure*}

        In Fig. \ref{fig:sim3:time_ev:v}, $V(x(k))$ of the constrained and unconstrained systems are compared.
        The maximal rest state $\bar{\bar x}$ might change after each modification, and so does $V(\bar{\bar x})$;  as expected, $V(x(k))\leq V(\bar x(k))$ in the unconstrained case. 
        Some modifications have longer adaption times and clearly the unconstrained system adapts faster than the constrained one; still, a global rest state $\bar x$ is reached faster than starting from an empty state, like at time $k=0$. 
        Removing nodes/arcs/sinks might increase the accumulation of tokens, as the length of the shortest paths might increase. 
        Conversely, inserting them might decrease the value of the state and  make it non-admissible (for a very brief time). 
        Removing sources has no effect, while inserting new ones makes their states reach the optimum values.

		\section{Conclusion}\label{sec:conclusion}
		We have considered the problem in which some agents are injected in an unknown network 
		and have to find a route to leave the network, moving using local information only. 
        We have proposed a decentralized local threshold
        policy, which is used by the tokens to decide whether to rest in the current node or move to some adjacent one.
        This mechanism has also been extended to support additional constraints on the paths that the tokens are allowed to follow.
        We have shown that, in the long run, tokens that apply our policy reach their aim and discover the shortest paths to leave, 
        even in the presence of an additional route constraint.
        We have also shown that the proposed policy is adaptive: when the network is dynamic, the shortest paths are updated if 
        the network maintains each new configuration  for a sufficiently long interval. 

  		\begin{IEEEbiography}[{\includegraphics[width=1in,height=1.25in,clip,keepaspectratio]{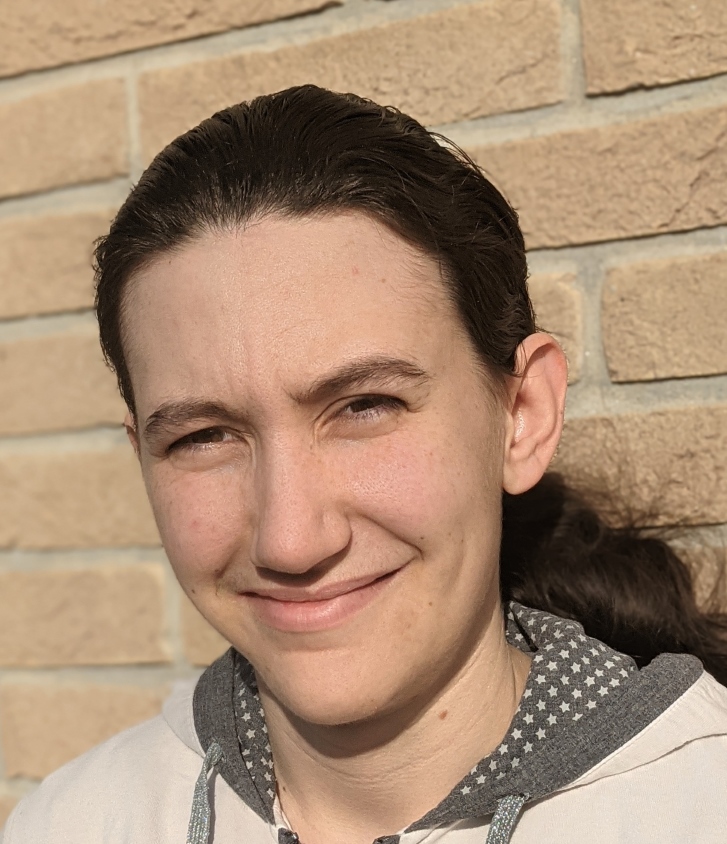}}]{Francesca Rosset}was born in S. Vito al Tagliamento (PN), Italy, in 1993. She received the B.Sc. and M.Sc. degree in electronic engineering in 2016 and 2019, respectively, both from the University of Udine, Italy. She is currently a PhD student  with the  Polytechnic Department of Engineering and Architecture of the University of Udine, Italy.
        Her interests concern the scheduling and optimization of smart grids and energy systems, and agent-based systems.
		\end{IEEEbiography}
		\vskip -2\baselineskip plus -1fil
		\begin{IEEEbiography}[{\includegraphics[width=1in,height=1.25in,clip,keepaspectratio]{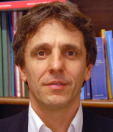}}]{Raffaele Pesenti}, 1964, is professor of Operations Research at the Department of Management of Universit\'a Ca' Foscari, Venice, Italy. His research interests concern the management and evaluation of complex systems with particular attention to the logistics and the biological ones. His papers deal with: management of distribution networks and systems, public transportation systems, game theory and consensus problem applications in production planning and control, logistics. He is Associate Editor for 4OR.
		\end{IEEEbiography}
		\vskip -2\baselineskip plus -1fil
        \begin{IEEEbiography}[{\includegraphics[width=1in,height=1.25in,clip,keepaspectratio]{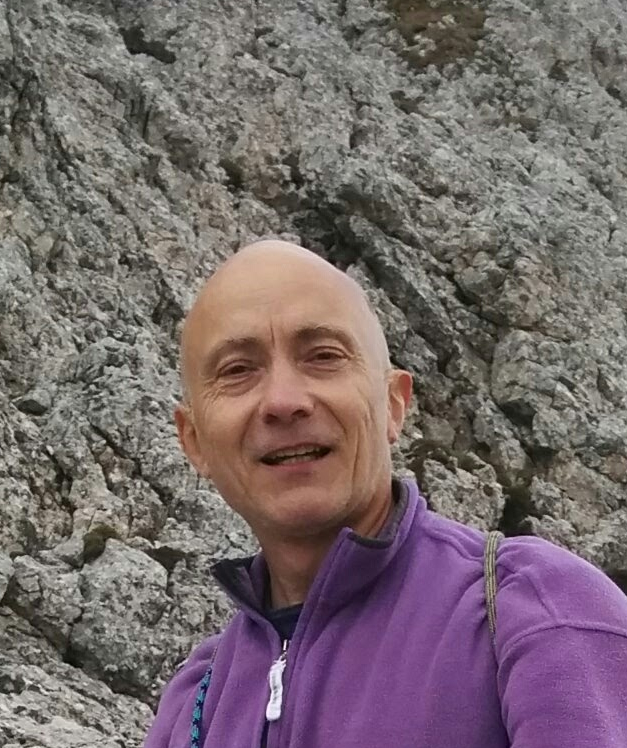}}]{Franco Blanchini}, 1959, is the Director of the Laboratory of System Dynamics at the University of Udine. He is co-author of the book ``Set theoretic methods in control'', Birkh\"auser. He received the 2001 ASME Oil \& Gas Application Committee Best Paper Award as a co--author of the article ``Experimental evaluation of a High-Gain Control for Compressor Surge Instability'', the 2002 IFAC prize Survey Paper Award as the author of the article ``Set Invariance in Control -- a survey'', Automatica, November 1999, for which he also received the High Impact Paper Award in 2017, and the 2017 NAHS Best Paper Award as a co-author of the article ``A switched system approach to dynamic race modelling'', Nonlinear Analysis: Hybrid Systems, 2016. He was nominated Senior Member of the IEEE in 2003. He was Associate Editor for Automatica and for IEEE Transactions on Automatic Control, and Senior Editor for IEEE Control Systems Letters.
		\end{IEEEbiography}
        \vfill

		\clearpage
		\appendix
        {\bf Proof of 	Lemma \ref{lem:adm_path}.} 
			\emph{If part}) If $x$ is admissible $x_{h_s}-x_{h_{s+1}} \leq \gamma_{h_s,h_{s+1}}$ for each pair  $h_s,h_{s+1} \in p$. Let $r$ be the number of (non-repeated) nodes in path $p$.
				Then,
                $$
				L(p)= \sum_{s=1}^{r-1}\gamma_{h_s,h_{s+1}}\geq \sum_{s=1}^{r-1}x_{h_s}-x_{h_{s+1}}=x_{h_1}-x_{h_r} = x_i-x_j,
				$$
                where the third term does not depend on the actual path, but just from the starting and the terminal node.
				
				\emph{Only if part}) If $x$ is not admissible there exists $(i,j) \in \mathcal{A}$ such that $x_{i}-x_{j} > \gamma_{ij}$.
                Then, for the trivial path $p =\{i,j\}$ we have $L(p) = \gamma_{ij} < x_{i}-x_{j}$, which completes the proof.\qed

        {\bf Proof of 	Lemma \ref{lem:Vbounded}.} 
			We first consider a network $\mathcal{G}$ with a single sink $j \in \mathcal{T}$,  maximum arc cost $\bar \gamma$, and state $x(k) \in \mathcal{O}$, for some $k \geq 0$. Recall that $x_j(k) \equiv 0$, by the definition of sink node.
			If  no token is injected, $v(k) = 0$, then  $x(k+1)= x(k)$.
			Also, for all paths~$p_i$ of $l_i$ arcs, joining any node $i \in \mathcal{N}$  with the sink $j$, it holds that $L(p_i)\leq l_i\bar\gamma$ and, by Lemma~\ref{lem:adm_path}, $x_i(k) \leq  L(p_i)$.   
			Now consider any anti-arborescence that joins all the nodes in~$\mathcal{N}$ to the sink $j$; for any $p_i$ in there, we have  $x_i(k) \leq l_i \bar{\gamma}$, and hence $V(x(k)) = \sum_{i \in \mathcal{N}}x_i(k) \leq \sum_{i \in \mathcal{N}} l_i \bar{\gamma}$.
			To complete the proof, note that for all the anti-arborescences it holds that: 
			$$\sum_{i \in \mathcal{N}} l_i \bar{\gamma} \leq \bar{\gamma}|\mathcal{N}|(|\mathcal{N}|-1)/2,$$ 
			where the equality holds if the anti-arborescence is a directed path. With multiple sinks, the proof is an easy but cumbersome generalization of this one: we have $x_j=0$ for all $j\in\mathcal{T}$, and the above argument holds for each anti-arborescence joining all the nodes to each sink node.\qed

        {\bf Proof of 	Theorem \ref{th:proc_defined}.} 
			To prove the theorem, we show that if at a time $k$ we have $x(k)\in \mathcal{O}$,
			then we get $x(k+1)\in\mathcal{O}$.
			
			If no token is injected at $k$, then  $x(k+1) = x(k) \in \mathcal{O}$.
   
			Instead, if a token is injected in a source $i \in \mathcal{S}$, we note that Policy \ref{pol:threshold} imposes that
			only the injected token may move through the network. Hence, two cases may occur.
   
		\emph{Case $x(k) + e_i \in \mathcal{O}$).} 
            By \eqref{eq:threshold_1}   and \eqref{eq:threshold_u},
			 we have that state equation~\eqref{eq:stateEqNode} 
			imposes $x(k+1) = x(k) + e_i \in \mathcal{O}$.
			
		\emph{Case $x(k) + e_i \not\in \mathcal{O}$).}
			Assume by contradiction that a state $x(k+1)\in \mathcal{O}$ is not eventually reached.
			This means that the token makes a walk of an infinite number of elementary transitions.
			Since the number of nodes is finite, such walk must include at least a circuit
			$p = \{i_1, i_2,\dots,i_r\equiv i_1\}$.
			As $u_{i_{h}i_{h+1}}=1$, for $h = 1,\ldots, r-1$, it implies  
			$x_{i_{h}} -x_{i_{h+1}} \geq \gamma_{i_{h},i_{h+1}}.$ 
			Summing up, we obtain
				\begin{equation*}
							\sum_{h=1}^{r-1}~x_{i_{h}} -x_{i_{h+1}} = x_{i_1}-x_{i_r}=  x_{i_1}-x_{i_1} = 0 \geq \sum_{h=1}^{r-1}\gamma_{i_{h},i_{h+1}}.
				\end{equation*}
			Hence, the length  $L(p)$ of the circuit $p$ is non-positive,
			in contradiction with Assumption~\ref{ass:nonNegCircuits}.2 on absence of non-positive circuits.
			As the injected token walk has no circuit, it is a path from the source $i$ to the token destination node $j$.
			Hence, $x(k+1) = x(k) \in \mathcal{O}$ if $j$ is a sink node, otherwise
			$x(k+1) = x(k) + e_j$. 
			In the latter case $x(k+1) \in \mathcal{O}$ since $x_r(k+1) = x_r(k)$ for all nodes $r \not = j$ and because
			Policy \ref{pol:threshold} imposes that if the token stops in $j$ then $x_j(k) + 1 - x_r(k) = x_j(k+1) - x_r(k)  \leq \gamma_{jk}$ for all $r \in \mathcal{N}_j$.\qed
			
		 {\bf Proof of 	Theorem \ref{th:proc_defined_2}.} 	
			From Theorem~\ref{th:proc_defined}, 
			the number of tokens in the network is non--decreasing, i.e, $V(x(k+1)) \geq V(x(k)$). 
   By Lemma~\ref{lem:Vbounded}, $V(x(k))$ converges to a finite value $\bar{V}$ from {below} and cannot diverge to $+\infty$.
			Also, as $V(x(k))$ is integer valued, $\bar{V}$ is reached in a finite time $\bar k$:
			\begin{equation}\label{eq:stadystate}
				 V(x(k)) = \bar{V} \leq  \bar{\gamma}|\mathcal{N}|(|\mathcal{N}|-1)/2,~~~k \geq \bar k.
        \vspace*{-2mm}
			\end{equation}
   
			If at some time $k'$ we reach a condition where injecting a new token in a {\em non-empty subset} of $\mathcal{S}$ always triggers the expulsion of a token from a sink, 
            by Theorem~\ref{th:proc_defined}, $x(k'+1)=x(k')$, and $V(x(k'+1))=V(x(k'))$; then, $x(\bar{k})$ is a rest state in $\partial\mathcal{O}$. By \eqref{eq:stadystate}, $V(x(k'))\leq \bar V$ and $k'\leq \bar k$ is finite.
          
		    Now assume that 
            $V(x(k))$ has converged to $\bar V$ by persistently injecting at least a token in {\em each} source.
		    Then, for $k \geq \bar{k}$, the injection of a new token in {\em any} source triggers 
		    the expulsion of a token from a sink, so $x(k+1)=x(k)=x(\bar{k})$, which means that $x(\bar{k})$ is a global rest state.
		    
		    In a single-source network, initially, every injected token increases $V(x(k))$ by $1$, 
            until $V(x(k))=\bar V$.
      Then, the maximum number of tokens to be inserted is  equal to the upper bound of $\bar V$, minus $V(x(0))$, the  number of tokens that are initially already present in the nodes' buffers. \qed

         {\bf Proof of 	Theorem \ref{th:proc_defined_2a}.}
			Under the theorem hypotheses, $V(x(k))$ is non-increasing, 
			since tokens can only leave the network.
			Then, either $V(x(k))$ converges to $\underline{V}$ from above or to $-\infty$.
			The latter situation may not occur. 
			Indeed, Policy \ref{pol:threshold} imposes that a
			token in node $i \in \mathcal{N}$ may leave the network only  
			along a path $p = \{i=h_1, h_2,\dots,h_r\}$, with $h_r\in\mathcal{T}$, 
			such that conditions $x_{h_1}-x_{h_2} > \gamma_{{h_1},h_{2}}$ and  $x_{h_s}-x_{h_{s+1}} \geq \gamma_{{h_s},h_{s+1}}$ for $s = 2,\dots r-1$, hold.
			Summing these up like in the proof of Theorem $\ref{th:proc_defined}$, as $x_{h_r}=0$ for sink $h_r$, we get that
			$x_i > L(p)$ must hold for a token in $i$ to exit: 
		 as the number of nodes is finite, the value of $V(x(k))$ cannot keep on decreasing.
			Differently, $V(x(k))$  converges in finite time to $\underline{V}$ as $V(x(k))$ is integer valued.
			
			Following the same reasoning in the proof of Theorem~\ref{th:proc_defined}, the
			subnetwork of~$\mathcal{G}$ induced by the arcs~$(i,j)$ traversed by a moving token 
			consists of a forest of non-intersecting paths and includes no circuits. 
			Hence, if $V(x(\underline{k})) = \underline{V}$ for  $\underline{k} \geq 0$ then
			$x(k) =  x(\underline{k})$ for all $k > \underline{k}$ if no token is injected.
			Then,  $x(\underline{k}) \in \mathcal{O}$.\qed

\begingroup
\tabcolsep=0.04cm
\begin{table*}[htb!]
\centering
\begin{threeparttable}
\caption{(\textit{Example 3}). Comparison of the simulation results at global rest state, for different small-world networks.\vspace*{-3mm}}
\label{tab:compare}
\scriptsize
\begin{tabular}{@{}cr rrrrrr rrc rrrrrr rrr r@{}}
\toprule
\multicolumn{2}{l}{\multirow{ 2}{*}{\thead{Parameter\\to vary}}}     & \multicolumn{6}{c}{{\footnotesize Unconstrained, original}} & \multicolumn{3}{c}{{\footnotesize Unc.,  enhanced}} & \multicolumn{6}{c}{{\footnotesize Constrained, original}}& \multicolumn{3}{c}{{\footnotesize Constr., enhanced}} & \multirow{ 2}{*}{\thead{$C_{max}$}}\\
\cmidrule(lr){3-8}\cmidrule(lr){9-11}\cmidrule(lr){12-17}\cmidrule(lr){18-20}
                             &            &{$L_{ss}$}&{$C_{ss}$}&{$E_{ss}$}&{$T_{ss}$}&{$V_{ss}$}&{$l_{ss}$}    &{$T_{ss}$}&{$V_{ss}$}&{$l_{ss}$}&       {$L_{ss}$}&{$C_{ss}$}&{$E_{ss}$}&{$T_{ss}$}&{$V_{ss}$}&{$l_{ss}$}      &{$T_{ss}$}&{$V_{ss}$}&{$l_{ss}$} &  \\
\midrule
\multirow{4}{*}{$n$}         & 100        & 104      & 40       & 7        & 1361     & 1344     & 1344         & 17       & 4509     & 0        &        122      & 28       &  6       & 9010     & 8886     & 8886           & 1042     & 15582    & 1018      & 32\\
                             & 1000       & 167      & 95       & 14       & 16427    & 16245    & 16245        & 64       & 78522    & 0        &        175      & 64       & 10       & 281645   & 279459   & 279459         & 22717    & 494319   & 22382     & 65\\
                             & 10000      & 218      & 125      & 20       & 159774   & 158746   & 158746       & 248      &1075638   & 0        &        229      & 81       & 12       & 3024376  & 3013981  & 3013981        & 240217   & 5325649  & 239205    & 85\\
                             & 50000      & 277      & 199      & 31       & 710604   & 709164   & 709164       & 341      &6993467   & 0        &        302      & 123      & 19       &44519968  & 43769899 & 43769899       &2436421   &86177388  &2387745    & 136\\  
\midrule
\multirow{5}{*}{$\delta$}    & 4          & 167      & 95       & 14       & 16427    & 16245    & 16245        & 64       & 78522    & 0        &        167      & 95       & 14       & 315195   & 309832   & 309832         & 12163    & 978229   & 11901     & 100\\
                             & 8          &  69      & 52       & 7        &  6786    &  6632    &  6632        & 42       & 37282    & 0        &         69      & 52       & 10       &  42498   &  41621   &  41621         &  3819    & 359209   &  3688     & 100\\
                             & 16         &  31      & 62       & 10       &  3325    &  3283    &  3283        & 36       & 17653    & 0        &         31      & 62       & 7        &  19699   &  19389   &  19389         &  3370    & 190531   &  3289     & 100\\
                             & 32         &  16      & 54       & 8        &  1956    &  1895    &  1895        & 64       & 8488     & 0        &         16      & 54       & 8        &   7483   &   7301   &   7301         &  2494    &  98433   &  2370     & 100\\
\midrule
\multirow{4}{*}{$E_{sp}$}       & 5          & 94       & 31       & 5        & 1155     & 1141     & 1141         & 12       & 41419    & 0        &        94       & 31       & 5        & 6712     & 6552     & 6552           & 355      & 46851    & 332       & 100\\
                             & 11         & 130      & 74       & 11       & 4699     & 4616     & 4616         & 30       &  60187   & 0        &        130      & 74       & 11       & 58029    & 56692    & 56692          & 2470     & 258335   & 2326      & 100\\
                             & 15         & 201      & 102      & 15       & 39993    & 39671    & 39671        & 71       & 110753   & 0        &        206      & 89       & 13       & 1261169  & 1229403  & 1229403        & 47694    & 2455789  & 46391     & 100\\
                             & 19         & 259      & 113      & 19       & 95199    & 94795    & 94795        & 110      & 157699   & 0        &        264      & 100      & 17       & 4223367  & 4174069  & 4174069        & 170082   & 5514798  & 167605    & 100\\
\midrule
\multirow{4}{*}{$K_{\gamma}$}& 1          & 167      & 95       & 14       & 16427    & 16245    & 16245        & 64       & 78522    & 0        &        175      & 64       & 10       & 281645   & 279459   & 279459         & 22717    & 494319   & 22382     & 65\\
                             & 10         & 1670     & 95       & 14       & 158364   & 158148   & 158148       & 64       &785220    & 0        &        1750     & 64       & 10       & 2716498  & 2714274  & 2714274        & 172647   & 4613687  & 172217    & 65\\
                             & 100        & 16700    & 95       & 14       & 1577358  & 1577178  & 1577178      & 64       &7852200   & 0        &        17500    & 64       & 10       & 27064645 & 27062424 & 27062424       &1675252   &45967239  & 1674833   & 65\\
\midrule
\multirow{4}{*}{$K_\sigma$}  & 1          & -        & -        &   -      & -        & -        & -            & -        & -        & -        &        175      & 64       &    10    & 281645   & 279459   & 279459         & 22717    & 494319   & 22382     & 65\\
                             & 10         & -        & -        &  -       & -        & -        & -            & -        & -        & -        &        175      & 640      &    10    & 274952   & 272982   & 272982         & 22789    & 476523   & 22386     & 645\\
                             & 100        & -        & -        &   -      & -        & -        & -            & -        & -        & -        &        175      & 6400     &    10    & 274952   & 272982   & 272982         & 22789    & 476523   & 22386     & 6450\\
\midrule
\multirow{4}{*}{\tiny $C_{max}$}   & 34         & -        & -        &  -       & -        & -        & -            & -        & -        & -        &        286      & 34       & 10       & 274704   & 274324   & 274324         & 72998    & 291584   & 72813     & 34\\
                             & 65         & -        & -        & -        & -        & -        & -            & -        & -        & -        &        175      & 64       & 10       & 281645   & 279459   & 279459         & 22717    & 494319   & 22382     & 65\\
                             & 95         & -        & -        & -        & -        & -        & -            & -        & -        & -        &        167      & 95       & 14       & 306771   & 301805   & 301805         & 12381    & 889270   & 12138     & 95\\
                             & 120        & -        & -        &  -       & -        & -        & -            & -        & -        & -        &        167      & 95       & 14       & 339537   & 333303   & 333303         & 11441    & 1342708  & 11129     & 120\\
\bottomrule
\end{tabular}%
\begin{tablenotes}[flushleft, para]\scriptsize
\note{Each row refers to a different network. Default parameters not specified in each row: number of nodes $n=1000$, mean nodes out-degree $\delta=4$, rewiring probability $\beta=0.15$, maximum costs $\gamma_{max}=50$, $\sigma_{max}=10$, costs' scaling factors $K_{\gamma}=1$, $K_{\sigma}=1$. $E_{sp}$: number of arcs in the shortest paths. $C_{max}$ is given for the constrained system.
{\bf Metrics.}  $L_{ss}$: cost $L(p)$ paid by each token injected at global rest state;  $C_{ss}$: the corresponding secondary cost $C(p)$; $E_{ss}$: the corresponding number of traversed arcs; $T_{ss}$: the time in which the global rest state is reached, i.e. $k$ such that $V(x(k))=\bar V$; $V_{ss}$: the corresponding value of $V(x(k))$; $l_{ss}$:  number of tokens lost during the transitory,  because  they are stopped in a node or their path is not feasible. $L_{ss}$, $C_{ss}$, $E_{ss}$ for the enhanced unconstrained (resp. constrained) policy are equal to the ones for the original policy, hence omitted.
}
    \end{tablenotes}
\end{threeparttable} 
\end{table*}
\endgroup

        {\bf Proof of 	Theorem \ref{th:shortestPath}.} 
			Assume that, at time $k$, $x(k)\in\mathcal{O}$, and a new token is injected
			in node $i\in\mathcal{S}$ and reaches its destination node $j$ along path $p=\{i=h_{0}, h_{1}, \dots, h_{r}=j\}$.
			Like in the proof of Theorem \ref{th:proc_defined}, we get $x_i(k)-x_j(k) = \sum_{s=0}^{r-1}\gamma_{h_{s},h_{s+1}} = L(p)$.
			Now assume that an alternative path $p'=\{i=h^\prime_{0}, h^\prime_{1}, \dots, h^\prime_{r}=j\}$ shorter than $p$ exists, so that $L(p^\prime) < L(p) = x_i(k)-x_j(k).$
			Then, by Lemma \ref{lem:adm_path}, as \eqref{eq:path} does not hold for $p^\prime$, the state $x$ is not admissible in contradiction with the hypothesis $x(k) \in \mathcal{O}$.
   The second part of the theorem holds as $x_j(k)\equiv0$ for all sinks $j\in\mathcal{T}$.\qed
		
        {\bf Proof of 	Theorem \ref{th:ext_acyclic}.} 
			Suppose that any circuit $\varphi$ in $\mathcal{G}$ fulfills $C(\varphi)>0$. 
			A token with secondary cost $c$ in node $i\in \varphi$ traversing this circuit will reach node $i$ again with cost $c+C(\varphi)>c$: 
			then, in $\mathcal{G}_E$ the token, initially in node $i^c$, would reach a different node $i^{c+C(\varphi)}$, closer to $i^{C_{max}}$. So, any walk in $\mathcal{G}$ becomes a path in $\mathcal{G}_E$. 
			
			Now assume a circuit $\varphi$ in $\mathcal{G}$ has $C(\varphi)=0$ and, by Assumption~\ref{ass:nonNegCircuits}, $L(\varphi)$ > 0: the token will reach again node $i$ with the same cost $c$; then, in $\mathcal{G}_E$,  
			 it would reach the same node $i^{c}$. 
			So, any circuit with $C(\varphi)=0$ in $\mathcal{G}$ corresponds to at most $C_{max}$ circuits $\varphi_E$ in $\mathcal{G}_E$ with  $L(\varphi_E) = L(\varphi)>0$.
			
			Conversely, suppose there exists a circuit $\varphi_E$ in $\mathcal{G}_E$. 
			A token in node $i^c\in \varphi_E$ traveling this circuit will reach node $i^c$ again. 
			Then, it travels a circuit $\varphi$ in $\mathcal{G}$, too: starting at node $i$ with secondary cost $c$, it returns to $i$ with the same $c$; this is possible only if $C(\varphi)=0$. Moreover, $L(\varphi)=L(\varphi_E)>0$, by construction and by Assumption \ref{ass:nonNegCircuits}.\qed

        {\bf Proof of 	Theorem \ref{th:constrshortestpaths}.} 
			Let $p=\{h^{c_1}_1,\ldots,h^{c_r}_r\}$
			 be a shortest path in $\mathcal{G}_E$ from source $h^{c_1}_1\in{\mathcal{S}_E}$ to a sink.
			 By Theorem~\ref{th:shortestPath}, there cannot be another path $p'$ with
			 $L(p') < L(p)$ joining $h^{c_1}_1$ with any other sink.
			By construction, $p$ must correspond to a 
			feasible shortest path in~$\mathcal{G}$. 
   
   Indeed, assume, by contradiction, that $p$ corresponds to a walk $\omega$ in $\mathcal{G}$.    
   By Theorem~\ref{th:ext_acyclic}, 
			there must be two nodes $h^{c_{s'}}_{s'}, h^{c_{s''}}_{s''}\in{\mathcal{N}_E}$ in path $p=\{h^{c_1}_1,\dots,\allowbreak h^{c_{s'}}_{s'},\dots, \allowbreak h^{c_{s''}}_{s''},\dots, \allowbreak h^{c_r}_r\}$, with $1\leq s'<s''\leq r$, that are associated with the same node $h_{s'}=h_{s''}\in\mathcal{N}$ in $\mathcal{G}$, so that the subpath 
			$p_\varphi = \{h^{c_{s'}}_{s'}, \dots, \allowbreak h^{c_{s''}}_{s''}\}$ corresponds to a circuit $\omega_\varphi$ of  $\omega$,
			with $C(\omega_\varphi)>0$ and $L_\varphi\doteq L(\omega_\varphi)=L(p_\varphi)>0$ by Assumption~\ref{ass:nonNegCircuits}.
			Consider now path $p''=\{h^{c_1}_1,\dots,\allowbreak h^{c_{s'}}_{s'},\allowbreak h^{c_{s''+1}-L_\varphi}_{s''+1},\dots, \allowbreak h^{c_r-L_\varphi}_r\}$ of $\mathcal{G}_E$
			that in $\mathcal{G}$ corresponds to the walk $\omega'$ obtained by removing  circuit $\omega_\varphi$ from $\omega$.
			We have $L(p'')=L(p)-L_\varphi<L(p)$ and $C(p'')=C(p)-C(\omega_\varphi)<C(p)$. So $\omega'$ remains feasible for $\mathcal{G}$, 
			and $p''$ is a path of $\mathcal{G}_E$, shorter than $p$, which is not possible, by Theorem~\ref{th:ext_acyclic}. \qed

        \subsection*{Numerical evaluation and comparison for different networks}
        
We consider single-source-single-sink Watts-Strogatz small-world networks \cite{Wat98}, i.e., random graphs with properties like clustering and short average path length, with $n$ nodes, $\delta$  mean nodes out-degree, $n\delta$ directed arcs (for each arc $(i,j)$, we also include $(j,i)$), and positive  $\gamma_{ij},\sigma_{ij}\in\mathbb{N}$ assigned randomly.
We compare the performance varying the network features, for  both the original and enhanced policy from Remark \ref{rem:enhpolicy},  using  a stochastic choice model.  
       
 Table~\ref{tab:compare} reports the parameters and some metrics summarizing the results at global rest state. 
 The enhanced policy always finds the  optimal paths as the original policy (see metrics $L_{ss}$, $C_{ss}$ and $E_{ss}$), but in considerably less time, especially for the unconstrained system. While the number of tokens stopped in the nodes is larger, no injected token is lost in  unconstrained systems as those are all virtual; in constrained systems some injected tokens fall asleep (and are lost), but they are fewer than the tokens lost with the original policy (see metrics $T_{ss}, V_{ss}$ and $l_{ss}$).

 Increasing $n$ or $E_{sp}$ (the number of arcs in the shortest paths), or decreasing $\delta$, results in a larger converging time. For the original policy, increasing $\gamma_{ij}$ by scaling factor $K_{\gamma}$ worsens the performance; the unconstrained enhanced policy is unaffected, while this is not true in the constrained case as tokens might fall asleep. Finally, scaling both $\sigma_{ij}$ and $C_{max}$ by $K_{\sigma}$ does not vary the performance, as the corresponding expanded networks have the same topology, while increasing only $C_{max}$ degrades the performance for the original policy, but improves it for the enhanced one.

\end{document}